\documentclass[11pt]{amsart}
\usepackage{amsmath,amssymb}
\newtheorem*{theorem}{Theorem}
\newtheorem{proposition}{Proposition}
\newtheorem{corollary}[proposition]{Corollary}

\begin{document}

\title{Blow-up phenomena for the Yamabe equation}
\author{Simon Brendle}
\date{May 5, 2007}

\maketitle

\begin{abstract} 
Let $(M,g)$ be compact Riemannian manifold of dimension $n \geq 3$. A well-known conjecture states that the set of constant scalar curvature metrics in the conformal class of $g$ is compact unless $(M,g)$ is conformally equivalent to the round sphere. In this paper, we construct counterexamples to this conjecture in dimensions $n \geq 52$.
\end{abstract}

\section{Introduction}

Let $(M,g)$ be a compact Riemannian manifold of dimension $n \geq 3$. The Yamabe problem is concerned with finding metrics of constant scalar curvature in the conformal class of $g$. This problem can be reduced to a semi-linear elliptic PDE. Indeed, the metric $u^{\frac{4}{n-2}} \, g$ has constant scalar curvature $c$ if and only if 
\begin{equation} 
\label{yamabe.pde} 
\frac{4(n-1)}{n-2} \, \Delta_g u - R_g \, u + c \, u^{\frac{n+2}{n-2}} = 0, 
\end{equation}
where $\Delta_g$ is the Laplace operator with respect to $g$ and $R_g$ denotes the scalar curvature of $g$. Clearly, every solution of (\ref{yamabe.pde}) is a critical point of the functional 
\begin{equation} 
\label{yamabe.functional}
E_g(u) = \frac{\int_M \big ( \frac{4(n-1)}{n-2} \, |du|_g^2 + R_g \, u^2 \big ) \, dvol_g}{\big ( \int_M 
u^{\frac{2n}{n-2}} \, dvol_g \big )^{\frac{n-2}{n}}}. 
\end{equation}
It is well-known that the PDE (\ref{yamabe.pde}) has at least one positive solution for any choice of $(M,g)$. If $n \geq 6$ and $(M,g)$ is not locally conformally flat, this follows from results of T.~Aubin \cite{Aubin1}. The remaining cases were solved by R.~Schoen using the positive mass theorem \cite{Schoen1}. 

Solutions to (\ref{yamabe.pde}) are not usually unique. As an example, consider the product metric on $S^1(L) \times S^{n-1}(1)$. If $L$ is sufficiently small, then the Yamabe PDE has a unique solution. On the other hand, there are many non-minimizing solutions if $L$ is large. D.~Pollack \cite{Pollack} has used gluing techniques to construct high energy solutions on more general background manifolds: given any conformal class with positive Yamabe constant and any positive integer $N$, there exists a new 
conformal class which is close to the original one in the $C^0$-norm, and contains at least $N$ metrics of constant scalar curvature (see \cite{Pollack}, Theorem 0.1). 

It is an interesting question whether the set of all solutions to the Yamabe PDE is compact (in the $C^2$-topology, say). A well-known conjecture states that this should be true unless $(M,g)$ is conformally equivalent to the round sphere (see \cite{Schoen2},\cite{Schoen3},\cite{Schoen4}). 
This conjecture has been verified in low dimensions and in the locally conformally flat case: 
if $(M,g)$ is locally conformally flat, compactness follows from work of R.~Schoen \cite{Schoen2},\cite{Schoen3}. 
Moreover, Schoen proposed a strategy for proving the conjecture in the non-locally conformally flat case based on the Pohozaev identity. In \cite{Li-Zhu}, Y.Y.~Li and M.~Zhu \cite{Li-Zhu} followed this strategy to prove compactness in dimension $3$. O.~Druet \cite{Druet} proved the conjecture in dimensions $4$ and $5$. Recently, F.~Marques \cite{Marques} showed that compactness holds up to dimension $7$. The same result was obtained independently by Y.Y.~Li and L.~Zhang \cite{Li-Zhang}. Moreover, Li and Zhang showed that compactness holds in all dimensions provided that $|W_g(p)| + |\nabla W_g(p)| > 0$ for all $p \in M$. M.~Khuri, F.~Marques, and R.~Schoen \cite{Khuri-Marques-Schoen} proved compactness up to dimension $24$, assuming that the positive mass theorem holds.\footnote{T.~Aubin has recently claimed a general compactness theorem in all dimensions \cite{Aubin2},\cite{Aubin3}. We have, however, been unable to verify some of the arguments in \cite{Aubin2}.}

In this paper, we address the opposite question: is it possible to construct Riemannian manifolds $(M,g)$ such that the set 
of constant scalar curvature metrics in the conformal class of $g$ is non-compact? So far, the only known examples where compactness fails involve non-smooth background metrics. The first result in this direction was established by A.~Ambrosetti and A.~Malchiodi \cite{Ambrosetti-Malchiodi}. This result was subsequently improved by M.~Berti and A.~Malchiodi \cite{Berti-Malchiodi}. Given positive integers $n$ and $k$ such that $k \geq 2$ and $n \geq 4k+3$, Berti and Malchiodi showed that there exists a Riemannian metric $g$ on $S^n$ (of class $C^k$) for which the set of solutions to the Yamabe PDE (\ref{yamabe.pde}) fails to be compact (see \cite{Berti-Malchiodi}, Theorem 1.1). A survey of these results can be found in \cite{Ambrosetti}. Recently, O.~Druet and E.~Hebey \cite{Druet-Hebey1} showed that blow-up can occur for problems of the form $Lu + c \, u^{\frac{n+2}{n-2}} = 0$, where $L$ is
a lower order perturbation of the conformal Laplacian on $S^n$.

We improve the results of Berti and Malchiodi by showing that the set of solutions to the Yamabe PDE (\ref{yamabe.pde}) can fail to be compact even if the background metric $g$ is $C^\infty$ smooth. In the examples we construct, the blowing-up sequence develops a singularity consisting of exactly one bubble.

\begin{theorem}
\label{main.theorem}
Assume that $n \geq 52$. Then there exists a Riemannian metric $g$ on $S^n$ (of class $C^\infty$) and a sequence of positive functions $v_\nu \in C^\infty(S^n)$ ($\nu \in \mathbb{N}$) with the following properties: 
\begin{itemize}
\item[(i)] $g$ is not conformally flat 
\item[(ii)] $v_\nu$ is a solution of the Yamabe PDE (\ref{yamabe.pde}) for all $\nu \in \mathbb{N}$ 
\item[(iii)] $E_g(v_\nu) < Y(S^n)$ for all $\nu \in \mathbb{N}$, and $E_g(v_\nu) \to Y(S^n)$ as $\nu \to \infty$ 
\item[(iv)] $\sup_{S^n} v_\nu \to \infty$ as $\nu \to \infty$
\end{itemize}
(Here, $Y(S^n)$ denotes the Yamabe energy of the round metric on $S^n$.)
\end{theorem}

Let us sketch the main steps involved in the proof of Theorem \ref{main.theorem}. For convenience, we will work on $\mathbb{R}^n$ instead of $S^n$. Let $g$ be a smooth metric on $\mathbb{R}^n$ which agrees with the Euclidean metric outside a ball of radius $1$. We will assume throughout the paper that $\det g(x) = 1$ for all $x \in \mathbb{R}^n$, so that the volume form associated with $g$ agrees with the Euclidean volume form. 

Our goal is to construct solutions to the Yamabe PDE on $(\mathbb{R}^n,g)$. In Section 2, we show that this problem can be reduced to finding critical points of a certain function $\mathcal{F}_g(\xi,\varepsilon)$, where $\xi$ is a vector in $\mathbb{R}^n$ and $\varepsilon$ is a positive real number. This idea has been used by many authors (see, e.g., \cite{Ambrosetti-Malchiodi} or \cite{Berti-Malchiodi}). In Section 3, we show that the function $\mathcal{F}_g(\xi,\varepsilon)$ can be approximated by an auxiliary function $F(\xi,\varepsilon)$. In Section 4, we prove that the function $F(\xi,\varepsilon)$ has a critical point, which is a strict local minimum. Finally, in Section 5, we use a perturbation argument to construct critical points of the function $\mathcal{F}_g(\xi,\varepsilon)$. From this the main result follows. 

It is a pleasure to thank Professor Richard Schoen and Professor Fernando Marques for numerous discussions on this subject. This project was supported by the Alfred P. Sloan Foundation and by the National Science Foundation under grant DMS-0605223. \\

\section{Lyapunov-Schmidt reduction}

Let
\[\mathcal{E} = \bigg \{ w \in L^{\frac{2n}{n-2}}(\mathbb{R}^n) \cap W_{loc}^{1,2}(\mathbb{R}^n): \int_{\mathbb{R}^n} |dw|^2 < \infty \bigg \}.\] 
By Sobolev's inequality, there exists a constant $K$, depending only on $n$, such that 
\[\bigg ( \int_{\mathbb{R}^n} |w|^{\frac{2n}{n-2}} \bigg )^{\frac{n-2}{n}} \leq K \, \int_{\mathbb{R}^n} |dw|^2\] 
for all $w \in \mathcal{E}$. We define a norm on $\mathcal{E}$ by $\|w\|_{\mathcal{E}}^2 = \int_{\mathbb{R}^n} |dw|^2$. It is easy to see that $\mathcal{E}$, equipped with this norm, is complete. \\

Given any pair $(\xi,\varepsilon) \in \mathbb{R}^n \times (0,\infty)$, we define a function $u_{(\xi,\varepsilon)}: 
\mathbb{R}^n \to \mathbb{R}$ by 
\[u_{(\xi,\varepsilon)}(x) = \Big ( \frac{\varepsilon}{\varepsilon^2 + |x - \xi|^2} \Big )^{\frac{n-2}{2}}.\] 
The function $u_{(\xi,\varepsilon)}$ satisfies the elliptic PDE 
\[\Delta u_{(\xi,\varepsilon)} + n(n-2) \, u_{(\xi,\varepsilon)}^{\frac{n+2}{n-2}} = 0.\] 
It is well known that 
\[\int_{\mathbb{R}^n} u_{(\xi,\varepsilon)}^{\frac{2n}{n-2}} = \Big ( \frac{Y(S^n)}{4n(n-1)} \Big )^{\frac{n}{2}}\] 
for all $(\xi,\varepsilon) \in \mathbb{R}^n \times (0,\infty)$. We next define 
\[\varphi_{(\xi,\varepsilon,0)}(x) = 
\Big ( \frac{\varepsilon}{\varepsilon^2 + |x - \xi|^2} \Big )^{\frac{n+2}{2}} \, \frac{\varepsilon^2 - |x - \xi|^2}{\varepsilon^2 + |x - \xi|^2}\] 
and 
\[\varphi_{(\xi,\varepsilon,k)}(x) = 
\Big ( \frac{\varepsilon}{\varepsilon^2 + |x - \xi|^2} \Big )^{\frac{n+2}{2}} \, \frac{2\varepsilon \, (x_k - \xi_k)}{\varepsilon^2 + |x - \xi|^2}\] 
for $k = 1,\hdots,n$. It is easy to see that the norm $\|\varphi_{(\xi,\varepsilon,k)}\|_{L^{\frac{2n}{n+2}}(\mathbb{R}^n)}$ is constant in $\xi$ and $\varepsilon$. Finally, we define a closed subspace $\mathcal{E}_{(\xi,\varepsilon)} \subset \mathcal{E}$ by 
\[\mathcal{E}_{(\xi,\varepsilon)} = \bigg \{ w \in \mathcal{E}: \int_{\mathbb{R}^n} \varphi_{(\xi,\varepsilon,k)} \, w = 0 \quad \text{for $k = 0,1,\hdots,n$} \bigg \}.\] 
Clearly, $u_{(\xi,\varepsilon)} \in \mathcal{E}_{(\xi,\varepsilon)}$.

\begin{proposition}
\label{bound.for.error.term}
Consider a Riemannian metric on $\mathbb{R}^n$ of the form $g(x) = \exp(h(x))$, where $h(x)$ is a trace-free symmetric two-tensor
on $\mathbb{R}^n$ satisfying $|h(x)| + |\partial h(x)| + |\partial^2 h(x)| \leq \alpha \leq 1$ for all $x \in \mathbb{R}^n$ and $h(x) = 0$ for $|x| \geq 1$. There exists a constant $C$, depending only on $n$, such that
\[\Big \| \Delta_g u_{(\xi,\varepsilon)} - \frac{n-2}{4(n-1)} \, R_g \, u_{(\xi,\varepsilon)} + n(n-2) \, u_{(\xi,\varepsilon)}^{\frac{n+2}{n-2}} \Big \|_{L^{\frac{2n}{n+2}}(\mathbb{R}^n)} \leq C \, \alpha\] 
for all pairs $(\xi,\varepsilon) \in \mathbb{R}^n \times (0,\infty)$.
\end{proposition}

\textbf{Proof.} 
Using the pointwise estimate 
\begin{align*} 
&\Big | \Delta_g u_{(\xi,\varepsilon)} - \frac{n-2}{4(n-1)} \, R_g \, u_{(\xi,\varepsilon)} + n(n-2) \, u_{(\xi,\varepsilon)}^{\frac{n+2}{n-2}} \Big | \\ 
&\leq C \, |h| \, |\partial^2 u_{(\xi,\varepsilon)}| + C \, |\partial h| \, |\partial u_{(\xi,\varepsilon)}| 
+ C \, (|\partial^2 h| + |\partial h|^2) \, u_{(\xi,\varepsilon)}, 
\end{align*} 
we obtain 
\begin{align*} 
&\Big \| \Delta_g u_{(\xi,\varepsilon)} - \frac{n-2}{4(n-1)} \, R_g \, u_{(\xi,\varepsilon)} + n(n-2) \, u_{(\xi,\varepsilon)}^{\frac{n+2}{n-2}} \Big \|_{L^{\frac{2n}{n+2}}(\mathbb{R}^n)} \\ 
&\leq C \, \|h\|_{L^\infty(\mathbb{R}^n)} \, \|\partial^2 u_{(\xi,\varepsilon)}\|_{L^{\frac{2n}{n+2}}(\mathbb{R}^n)} 
+ C \, \|\partial h\|_{L^n(\mathbb{R}^n)} \, \|\partial u_{(\xi,\varepsilon)}\|_{L^2(\mathbb{R}^n)} \\ 
&+ C \, 
(\|\partial^2 h\|_{L^{\frac{n}{2}}(\mathbb{R}^n)} + \|\partial h\|_{L^n(\mathbb{R}^n)}^2) \, \|u_{(\xi,\varepsilon)}\|_{L^{\frac{2n}{n-2}}(\mathbb{R}^n)} \\ 
&\leq C \, \alpha. 
\end{align*} 
This proves the assertion. \\

\begin{proposition} 
\label{eigenvalue.estimate.1}
There exists a positive constant $\theta$, depending only on $n$, such that 
\begin{align*} 
&\int_{\mathbb{R}^n} \Big ( |dw|^2 - n(n+2) \, u_{(\xi,\varepsilon)}^{\frac{4}{n-2}} \, w^2 \Big ) \\ 
&\geq 2\theta \, \|w\|_{\mathcal{E}}^2 - \frac{16n^2}{\theta} \, 
\bigg ( \int_{\mathbb{R}^n} u_{(\xi,\varepsilon)}^{\frac{n+2}{n-2}} \, w \bigg )^2 
\end{align*}
for all $w \in \mathcal{E}_{(\xi,\varepsilon)}$.
\end{proposition}

Proposition \ref{eigenvalue.estimate.1} follows from an analysis of the eigenvalues of the Laplace operator on $S^n$. The details can be found in \cite{Rey}. \\

\begin{corollary} 
\label{eigenvalue.estimate.2}
Consider a Riemannian metric on $\mathbb{R}^n$ of the form $g(x) = \exp(h(x))$, where $h(x)$ is a trace-free symmetric two-tensor on $\mathbb{R}^n$ satisfying $h(x) = 0$ for $|x| \geq 1$. 
There exists a positive constant $\alpha_0 \leq 1$, depending only on $n$, with the following property: if $|h(x)| + |\partial h(x)| + |\partial^2 h(x)| \leq \alpha_0$ for all $x \in 
\mathbb{R}^n$, then we have 
\[\bigg ( \int_{\mathbb{R}^n} |w|^{\frac{2n}{n-2}} \bigg )^{\frac{n-2}{n}} \leq 2K \, \int_{\mathbb{R}^n} 
\Big ( |dw|_g^2 + \frac{n-2}{4(n-1)} \, R_g \, w^2 \Big )\] 
for all $w \in \mathcal{E}$ and 
\begin{align*} 
&\int_{\mathbb{R}^n} \Big ( |dw|_g^2 + \frac{n-2}{4(n-1)} \, R_g \, w^2 - n(n+2) \, u_{(\xi,\varepsilon)}^{\frac{4}{n-2}} \, w^2 \Big ) \\ 
&\geq \theta \, \|w\|_{\mathcal{E}}^2 - \frac{1}{\theta} \, \bigg ( \int_{\mathbb{R}^n} \Big ( \Delta_g u_{(\xi,\varepsilon)} - \frac{n-2}{4(n-1)} \, R_g \, u_{(\xi,\varepsilon)} + n(n+2) \, u_{(\xi,\varepsilon)}^{\frac{n+2}{n-2}} \Big ) \, w \bigg )^2 
\end{align*}
for all $w \in \mathcal{E}_{(\xi,\varepsilon)}$.
\end{corollary}

\textbf{Proof.} 
Using Proposition \ref{bound.for.error.term} and H\"older's inequality, we obtain 
\begin{align*} 
&\bigg | \int_{\mathbb{R}^n} \Big ( \Delta_g u_{(\xi,\varepsilon)} - \frac{n-2}{4(n-1)} \, R_g \, u_{(\xi,\varepsilon)} + n(n+2) \, u_{(\xi,\varepsilon)}^{\frac{n+2}{n-2}} \Big ) \, w \bigg | \\ 
&\geq 4n \, \bigg | \int_{\mathbb{R}^n} u_{(\xi,\varepsilon)}^{\frac{n+2}{n-2}} \, w \bigg | - C \, \alpha_0 \, \|w\|_{\mathcal{E}}. 
\end{align*} 
This implies 
\begin{align*} 
&\bigg ( \int_{\mathbb{R}^n} \Big ( \Delta_g u_{(\xi,\varepsilon)} - \frac{n-2}{4(n-1)} \, R_g \, u_{(\xi,\varepsilon)} + n(n+2) \, u_{(\xi,\varepsilon)}^{\frac{n+2}{n-2}} \Big ) \, w \bigg )^2 \\ 
&\geq 16n^2 \, \bigg ( \int_{\mathbb{R}^n} u_{(\xi,\varepsilon)}^{\frac{n+2}{n-2}} \, w \bigg )^2 - \theta^2 \, \|w\|_{\mathcal{E}}^2 
\end{align*} 
if $\alpha_0$ is sufficiently small. Hence, the assertion follows from Proposition \ref{eigenvalue.estimate.1}. \\

\begin{proposition} 
\label{linearized.operator}
Consider a Riemannian metric on $\mathbb{R}^n$ of the form $g(x) = \exp(h(x))$, where $h(x)$ is a trace-free symmetric two-tensor on $\mathbb{R}^n$ satisfying $|h(x)| + |\partial h(x)| + |\partial^2 h(x)| \leq \alpha_0$ for all $x \in \mathbb{R}^n$ and $h(x) = 0$ for $|x| \geq 1$. 
Given any pair $(\xi,\varepsilon) \in \mathbb{R}^n \times (0,\infty)$ and any function $f \in L^{\frac{2n}{n+2}}(\mathbb{R}^n)$, there exists a unique function $w \in \mathcal{E}_{(\xi,\varepsilon)}$ such that
\[\int_{\mathbb{R}^n} \Big ( \langle dw,d\psi \rangle_g 
+ \frac{n-2}{4(n-1)} \, R_g \, w \, \psi - n(n+2) \, u_{(\xi,\varepsilon)}^{\frac{4}{n-2}} \, w \, \psi \Big ) = \int_{\mathbb{R}^n} f \, \psi\] 
for all test functions $\psi \in \mathcal{E}_{(\xi,\varepsilon)}$. Moreover, we have $\|w\|_{\mathcal{E}} \leq C \, \|f\|_{L^{\frac{2n}{n+2}}(\mathbb{R}^n)}$, where $C$ is a constant that depends only on $n$.
\end{proposition}

\textbf{Proof.} 
Suppose that $w \in \mathcal{E}_{(\xi,\varepsilon)}$ and 
\[\int_{\mathbb{R}^n} \Big ( \langle dw,d\psi \rangle_g 
+ \frac{n-2}{4(n-1)} \, R_g \, w \, \psi - n(n+2) \, u_{(\xi,\varepsilon)}^{\frac{4}{n-2}} \, w \, \psi \Big ) = \int_{\mathbb{R}^n} f \, \psi\] 
for all test functions $\psi \in \mathcal{E}_{(\xi,\varepsilon)}$. This implies 
\[\int_{\mathbb{R}^n} \Big ( |dw|_g^2 + \frac{n-2}{4(n-1)} \, R_g \, w^2 - n(n+2) \, u_{(\xi,\varepsilon)}^{\frac{4}{n-2}} \, w^2 \Big ) = \int_{\mathbb{R}^n} f \, w\] 
and 
\[\int_{\mathbb{R}^n} \Big ( \Delta_g u_{(\xi,\varepsilon)} - \frac{n-2}{4(n-1)} \, R_g \, u_{(\xi,\varepsilon)} + n(n+2) \, u_{(\xi,\varepsilon)}^{\frac{n+2}{n-2}} \Big ) \, w 
= -\int_{\mathbb{R}^n} u_{(\xi,\varepsilon)} \, f.\] 
Using Corollary \ref{eigenvalue.estimate.2}, we obtain 
\begin{align*} 
\theta \, \|w\|_{\mathcal{E}}^2 &\leq \int_{\mathbb{R}^n} \Big ( |dw|_g^2 + \frac{n-2}{4(n-1)} \, R_g \, w^2 - n(n+2) \, u_{(\xi,\varepsilon)}^{\frac{4}{n-2}} \, w^2 \Big ) \\ 
&+ \frac{1}{\theta} \, \bigg ( \int_{\mathbb{R}^n} \Big ( \Delta_g u_{(\xi,\varepsilon)} - \frac{n-2}{4(n-1)} \, R_g \, u_{(\xi,\varepsilon)} + n(n+2) \, u_{(\xi,\varepsilon)}^{\frac{n+2}{n-2}} \Big ) \, w \bigg )^2 \\ 
&\leq \bigg ( \int_{\mathbb{R}^n} |f|^{\frac{2n}{n+2}} \bigg )^{\frac{n+2}{2n}} \, \bigg ( \int_{\mathbb{R}^n} |w|^{\frac{2n}{n-2}} \bigg )^{\frac{n-2}{2n}} \\ 
&+ \frac{1}{\theta} \, \bigg ( \int_{\mathbb{R}^n} u_{(\xi,\varepsilon)}^{\frac{2n}{n-2}} \bigg )^{\frac{n-2}{n}} \, \bigg ( \int_{\mathbb{R}^n} |f|^{\frac{2n}{n+2}} \bigg )^{\frac{n+2}{n}} \\ 
&\leq K^{\frac{1}{2}} \, \|f\|_{L^{\frac{2n}{n+2}}(\mathbb{R}^n)} \, \|w\|_{\mathcal{E}} + \frac{1}{\theta} \, \Big ( \frac{Y(S^n)}{4n(n-1)} \Big )^{\frac{n-2}
{2}} \, \|f\|_{L^{\frac{2n}{n+2}}(\mathbb{R}^n)}^2.
\end{align*} 
Hence, it follows from Young's inequality that 
\[\frac{\theta}{2} \, \|w\|_{\mathcal{E}}^2 
\leq \frac{K}{2\theta} \, \|f\|_{L^{\frac{2n}{n+2}}(\mathbb{R}^n)}^2 + \frac{1}{\theta} \, \Big ( \frac{Y(S^n)}{4n(n-1)} \Big )^{\frac{n-2}{2}} \, \|f\|_{L^{\frac{2n}{n+2}}(\mathbb{R}^n)}^2.\] 
From this the uniqueness statement follows easily. 

In order to prove the existence part, it suffices to minimize the functional 
\begin{align*} 
&\int_{\mathbb{R}^n} \Big ( |dw|_g^2 + \frac{n-2}{4(n-1)} \, R_g \, w^2 - n(n+2) \, u_{(\xi,\varepsilon)}^{\frac{4}{n-2}} \, w^2 - 2fw \Big ) \\ 
&+ \frac{1}{\theta} \, \bigg ( \int_{\mathbb{R}^n} \Big ( \Delta_g u_{(\xi,\varepsilon)} - \frac{n-2}{4(n-1)} \, R_g \, u_{(\xi,\varepsilon)} + n(n+2) \, u_{(\xi,\varepsilon)}^{\frac{n+2}{n-2}} \Big ) \, w \bigg )^2 
\end{align*}
over all functions $w \in \mathcal{E}_{(\xi,\varepsilon)}$. \\

\begin{proposition} 
\label{fixed.point.argument}
Consider a Riemannian metric on $\mathbb{R}^n$ of the form $g(x) = \exp(h(x))$, where $h(x)$ is a trace-free symmetric two-tensor on $\mathbb{R}^n$ satisfying $h(x) = 0$ for $|x| \geq 1$. Moreover, let $(\xi,\varepsilon) \in \mathbb{R}^n \times (0,\infty)$. There exists a positive constant $\alpha_1 \leq \alpha_0$, depending only on $n$, with the following property: if $|h(x)| + |\partial h(x)| + |\partial^2 h(x)| \leq \alpha_1$ for all $x \in \mathbb{R}^n$, then there exists a function $v_{(\xi,\varepsilon)} \in \mathcal{E}$ such that $v_{(\xi,\varepsilon)} - u_{(\xi,\varepsilon)} \in \mathcal{E}_{(\xi,\varepsilon)}$ 
and 
\[\int_{\mathbb{R}^n} \Big ( \langle dv_{(\xi,\varepsilon)},d\psi \rangle_g + \frac{n-2}{4(n-1)} \, R_g \, v_{(\xi,\varepsilon)} \, \psi - n(n-2) \, |v_{(\xi,\varepsilon)}|^{\frac{4}{n-2}} \, v_{(\xi,\varepsilon)} \, \psi \Big ) = 0\] for all test functions $\psi \in \mathcal{E}_{(\xi,\varepsilon)}$. Moreover, we have the estimate 
\begin{align*} 
&\|v_{(\xi,\varepsilon)} - u_{(\xi,\varepsilon)}\|_{\mathcal{E}} \\ 
&\leq C \, \Big \| \Delta_g u_{(\xi,\varepsilon)} - \frac{n-2}{4(n-1)} \, R_g \, u_{(\xi,\varepsilon)} + n(n-2) \, u_{(\xi,\varepsilon)}^{\frac{n+2}{n-2}} \Big \|_{L^{\frac{2n}{n+2}}(\mathbb{R}^n)}, 
\end{align*}
where $C$ is a constant that depends only on $n$.
\end{proposition}

\textbf{Proof.} 
Let $G_{(\xi,\varepsilon)}: L^{\frac{2n}{n+2}}(\mathbb{R}^n) \to \mathcal{E}_{(\xi,\varepsilon)}$ be the solution operator constructed in Proposition \ref{linearized.operator}. We define a nonlinear mapping $\Phi_{(\xi,\varepsilon)}: \mathcal{E}_{(\xi,\varepsilon)} \to \mathcal{E}_{(\xi,\varepsilon)}$ by 
\begin{align*} 
&\Phi_{(\xi,\varepsilon)}(w) \\ 
&= G_{(\xi,\varepsilon)} \Big ( \Delta_g u_{(\xi,\varepsilon)} - \frac{n-2}{4(n-1)} \, R_g \, u_{(\xi,\varepsilon)} + n(n-2) \, u_{(\xi,\varepsilon)}^{\frac{n+2}{n-2}} \Big ) \\ 
&+ n(n-2) \, G_{(\xi,\varepsilon)} 
\Big ( |u_{(\xi,\varepsilon)} + w|^{\frac{4}{n-2}} \, (u_{(\xi,\varepsilon)} + w) - u_{(\xi,\varepsilon)}^{\frac{n+2}{n-2}} - \frac{n+2}{n-2} \, u_{(\xi,\varepsilon)}^{\frac{4}{n-2}} \, w \Big ). 
\end{align*} 
It follows from Proposition \ref{bound.for.error.term} that 
$\|\Phi_{(\xi,\varepsilon)}(0)\|_{\mathcal{E}} \leq C \, \alpha_1$. 
Using the pointwise estimate 
\begin{align*} 
&\Big | |u_{(\xi,\varepsilon)} + w|^{\frac{4}{n-2}} \, (u_{(\xi,\varepsilon)} + w) - |u_{(\xi,\varepsilon)} + \tilde{w}|^{\frac{4}{n-2}} \, (u_{(\xi,\varepsilon)} + \tilde{w}) \\ 
&\hspace{20mm} - \frac{n+2}{n-2} \, u_{(\xi,\varepsilon)}^{\frac{4}{n-2}} \, (w - \tilde{w}) \Big | \\ 
&\leq C \, (|w|^{\frac{4}{n-2}} + |\tilde{w}|^{\frac{4}{n-2}}) \, |w - \tilde{w}|, 
\end{align*} 
we obtain 
\begin{align*} 
&\|\Phi_{(\xi,\varepsilon)}(w) - \Phi_{(\xi,\varepsilon)}(\tilde{w})\|_{\mathcal{E}} \\ 
&\leq C \, \Big \| |u_{(\xi,\varepsilon)} + w|^{\frac{4}{n-2}} \, (u_{(\xi,\varepsilon)} + w) - |u_{(\xi,\varepsilon)} + \tilde{w}|^{\frac{4}{n-2}} \, (u_{(\xi,\varepsilon)} + \tilde{w}) \\ 
&\hspace{20mm} - \frac{n+2}{n-2} \, u_{(\xi,\varepsilon)}^{\frac{4}{n-2}} \, (w - \tilde{w}) \Big \|_{L^{\frac{2n}{n+2}}(\mathbb{R}^n)} \\ 
&\leq C \, (\|w\|_{L^{\frac{2n}{n-2}}(\mathbb{R}^n)}^{\frac{4}{n-2}} + \|\tilde{w}\|_{L^{\frac{2n}{n-2}}(\mathbb{R}^n)}^{\frac{4}{n-2}}) \, \|w - \tilde{w}\|_{L^{\frac{2n}{n-2}}(\mathbb{R}^n)} 
\end{align*} 
for all functions $w,\tilde{w} \in \mathcal{E}_{(\xi,\varepsilon)}$. This implies 
\[\|\Phi_{(\xi,\varepsilon)}(w) - \Phi_{(\xi,\varepsilon)}(\tilde{w})\|_{\mathcal{E}} \leq C \, (\|w\|_{\mathcal{E}}^{\frac{4}{n-2}} + \|\tilde{w}\|_{\mathcal{E}}^{\frac{4}{n-2}}) \, \|w - \tilde{w}\|_{\mathcal{E}}\] 
for $w,\tilde{w} \in \mathcal{E}_{(\xi,\varepsilon)}$. Hence, if $\alpha_1$ is sufficiently small, then the contraction mapping principle implies that the mapping $\Phi_{(\xi,\varepsilon)}$ has a unique fixed point. From this the assertion follows easily. \\

We next define a function $\mathcal{F}_g: \mathbb{R}^n \times (0,\infty) \to \mathbb{R}$ by 
\begin{align*} 
\mathcal{F}_g(\xi,\varepsilon) &= \int_{\mathbb{R}^n} 
\Big ( |dv_{(\xi,\varepsilon)}|_g^2 + \frac{n-2}{4(n-1)} \, R_g \, v_{(\xi,\varepsilon)}^2 - (n-2)^2 \, |v_{(\xi,\varepsilon)}|^{\frac{2n}{n-2}} \Big ) \\ 
&- 2(n-2) \, \Big ( \frac{Y(S^n)}{4n(n-1)} \Big )^{\frac{n}{2}}. 
\end{align*} 
If we choose $\alpha_1$ small enough, then we obtain the following result: \\

\begin{proposition} 
\label{reduction.to.a.finite.dimensional.problem}
The function $\mathcal{F}_g$ is continuously differentiable. Moreover, if $(\bar{\xi},\bar{\varepsilon})$ is a critical point of the function $\mathcal{F}_g$, then 
the function $v_{(\bar{\xi},\bar{\varepsilon})}$ is a non-negative weak solution of the equation 
\[\Delta_g v_{(\bar{\xi},\bar{\varepsilon})} - \frac{n-2}{4(n-1)} \, R_g \, v_{(\bar{\xi},\bar{\varepsilon})} + n(n-2) \, v_{(\bar{\xi},\bar{\varepsilon})}^{\frac{n+2}{n-2}} = 0.\] 
\end{proposition}

\textbf{Proof.} 
By definition of $v_{(\xi,\varepsilon)}$, we can find real numbers $a_k(\xi,\varepsilon)$, $k = 0,1, \hdots, n$, such that 
\begin{align*} 
&\int_{\mathbb{R}^n} \Big ( \langle dv_{(\xi,\varepsilon)},d\psi \rangle_g + \frac{n-2}{4(n-1)} \, R_g \, v_{(\xi,\varepsilon)} \, v_{(\xi,\varepsilon)} \, \psi - n(n-2) \, |v_{(\xi,\varepsilon)}|^{\frac{4}{n-2}} \, v_{(\xi,\varepsilon)} \, \psi \Big ) \\ 
&= \sum_{k=0}^n a_k(\xi,\varepsilon) \, \int_{\mathbb{R}^n} \varphi_{(\xi,\varepsilon,k)} \, \psi 
\end{align*} 
for all test functions $\psi \in \mathcal{E}$. This implies 
\[\frac{\partial}{\partial \varepsilon} \mathcal{F}_g(\varepsilon,\xi) = 2 \, \sum_{k=0}^n a_k(\xi,\varepsilon) \, \int_{\mathbb{R}^n} \varphi_{\xi,\varepsilon,k)} \, \frac{\partial}{\partial \varepsilon} v_{(\xi,\varepsilon)}\] 
and 
\[\frac{\partial}{\partial \xi_j} \mathcal{F}_g(\varepsilon,\xi) = 2 \, \sum_{k=0}^n a_k(\xi,\varepsilon) \, \int_{\mathbb{R}^n} \varphi_{(\xi,\varepsilon,k)} \, \frac{\partial}{\partial \xi_j} v_{(\xi,\varepsilon)}\] 
for $j = 1, \hdots, n$. On the other hand, we have 
\[\int_{\mathbb{R}^n} \varphi_{(\xi,\varepsilon,k)} \, (v_{(\xi,\varepsilon)} - u_{(\xi,\varepsilon)}) = 0\] 
since $v_{(\xi,\varepsilon)} - u_{(\xi,\varepsilon)} \in \mathcal{E}_{(\xi,\varepsilon)}$. This implies 
\begin{align*} 
0 &= \int_{\mathbb{R}^n} \frac{\partial}{\partial \varepsilon} \varphi_{(\xi,\varepsilon,k)} \, (v_{(\xi,\varepsilon)} - u_{(\xi,\varepsilon)}) + \int_{\mathbb{R}^n} \varphi_{(\xi,\varepsilon,k)} \, \frac{\partial}{\partial \varepsilon} (v_{(\xi,\varepsilon)} - u_{(\xi,\varepsilon)}) \\ 
&= \int_{\mathbb{R}^n} \frac{\partial}{\partial \varepsilon} \varphi_{(\xi,\varepsilon,k)} \, (v_{(\xi,\varepsilon)} - u_{(\xi,\varepsilon)}) + \int_{\mathbb{R}^n} \varphi_{(\xi,\varepsilon,k)} \, \frac{\partial}{\partial \varepsilon} v_{(\xi,\varepsilon)} \\ 
&+ \frac{n-2}{2(n+1)} \, \Big ( \frac{Y(S^n)}{4n(n-1)} \Big )^{\frac{n}{2}} \, \varepsilon^{-1} \, \delta_{0k} 
\end{align*} 
and 
\begin{align*} 
0 &= \int_{\mathbb{R}^n} \frac{\partial}{\partial \xi_j} \varphi_{(\xi,\varepsilon,k)} \, (v_{(\xi,\varepsilon)} - u_{(\xi,\varepsilon)}) + \int_{\mathbb{R}^n} \varphi_{(\xi,\varepsilon,k)} \, \frac{\partial}{\partial \xi_j} (v_{(\xi,\varepsilon)} - u_{(\xi,\varepsilon)}) \\ 
&= \int_{\mathbb{R}^n} \frac{\partial}{\partial \xi_j} \varphi_{(\xi,\varepsilon,k)} \, (v_{(\xi,\varepsilon)} - u_{(\xi,\varepsilon)}) + \int_{\mathbb{R}^n} \varphi_{(\xi,\varepsilon,k)} \, \frac{\partial}{\partial \xi_j} v_{(\xi,\varepsilon)} \\ 
&- \frac{n-2}{2(n+1)} \, \Big ( \frac{Y(S^n)}{4n(n-1)} \Big )^{\frac{n}{2}} \, \varepsilon^{-1} \, \delta_{jk} 
\end{align*} 
for $j = 1,\hdots,n$. Putting these facts together, we obtain 
\begin{align*} 
&- \frac{n-2}{n+1} \, \Big ( \frac{Y(S^n)}{4n(n-1)} \Big )^{\frac{n}{2}} \, a_0(\xi,\varepsilon) \\ 
&= \varepsilon \, \frac{\partial}{\partial \varepsilon} \mathcal{F}_g(\xi,\varepsilon) + 2\varepsilon \, \sum_{k=0}^n a_k(\xi,\varepsilon) \, \int_{\mathbb{R}^n} \frac{\partial}{\partial \varepsilon} \varphi_{(\xi,\varepsilon,k)} \, (v_{(\xi,\varepsilon)} - u_{(\xi,\varepsilon)}) 
\end{align*} 
and 
\begin{align*} 
&\frac{n-2}{n+1} \, \Big ( \frac{Y(S^n)}{4n(n-1)} \Big )^{\frac{n}{2}} \, a_j(\xi,\varepsilon) \\ 
&= \varepsilon \, \frac{\partial}{\partial \xi_j} \mathcal{F}_g(\xi,\varepsilon) + 2\varepsilon \, \sum_{k=0}^n a_k(\xi,\varepsilon) \, \int_{\mathbb{R}^n} \frac{\partial}{\partial \xi_j} \varphi_{(\xi,\varepsilon,k)} \, (v_{(\xi,\varepsilon)} - u_{(\xi,\varepsilon)}) 
\end{align*} 
for $j = 1, \hdots, n$. Hence, if $(\bar{\xi},\bar{\varepsilon})$ is a critical point of $\mathcal{F}_g$, then we have 
\[\sum_{k=0}^n |a_k(\bar{\xi},\bar{\varepsilon})| \leq C \, \|v_{(\bar{\xi},\bar{\varepsilon})} - u_{(\bar{\xi},\bar{\varepsilon})}\|_{L^{\frac{2n}{n-2}}(\mathbb{R}^n)} \, \sum_{k=0}^n |a_k(\bar{\xi},\bar{\varepsilon})|,\] 
where $C$ is a constant that depends only on $n$.
On the other hand, we have $\|v_{(\bar{\xi},\bar{\varepsilon})} - u_{(\bar{\xi},\bar{\varepsilon})}\|_{L^{\frac{2n}{n-2}}(\mathbb{R}^n)} \leq C \, \alpha_1$.
Hence, if we choose $\alpha_1$ sufficiently small, then we must have $a_k(\bar{\xi},\bar{\varepsilon}) = 0$ for $k = 0,1, \hdots, n$. Thus, we conclude that
\[\int_{\mathbb{R}^n} \Big ( \langle dv_{(\bar{\xi},\bar{\varepsilon})},d\psi \rangle_g + \frac{n-2}{4(n-1)} \, R_g \, v_{(\bar{\xi},\bar{\varepsilon})} \, \psi - n(n-2) \, |v_{(\bar{\xi},\bar{\varepsilon})}|^{\frac{4}{n-2}} \, v_{(\bar{\xi},\bar{\varepsilon})} \, \psi \Big ) = 0\] 
for all test functions $\psi \in \mathcal{E}$. 
It remains to show that the function $v_{(\bar{\xi},\bar{\varepsilon})}$ is non-negative. 
To that end, we put $\psi = \min \{v_{(\bar{\xi},\bar{\varepsilon})},0\}$. Since 
$v_{(\bar{\xi},\bar{\varepsilon})} \in \mathcal{E}$, we conclude that $\psi \in \mathcal{E}$. 
This implies 
\begin{align*} 
&\int_{\{v_{(\bar{\xi},\bar{\varepsilon})} < 0\}} \Big ( |dv_{(\bar{\xi},\bar{\varepsilon})}|_g^2 + \frac{n-2}{4(n-1)} \, R_g \, v_{(\bar{\xi},\bar{\varepsilon})}^2 \Big ) \\ 
&= n(n-2) \, \int_{\{v_{(\bar{\xi},\bar{\varepsilon})} < 0\}} |v_{(\bar{\xi},\bar{\varepsilon})}|^{\frac{2n}{n-2}}. 
\end{align*} 
Moreover, we have 
\begin{align*} 
&\bigg ( \int_{\{v_{(\bar{\xi},\bar{\varepsilon})} < 0\}} |v_{(\bar{\xi},\bar{\varepsilon})}|^{\frac{2n}{n-2}} \bigg )^{\frac{n-2}{n}} \\ 
&\leq 2K \int_{\{v_{(\bar{\xi},\bar{\varepsilon})} < 0\}} \Big ( |dv_{(\bar{\xi},\bar{\varepsilon})}|_g^2 + \frac{n-2}{4(n-1)} \, R_g \, v_{(\bar{\xi},\bar{\varepsilon})}^2 \Big ) 
\end{align*} 
by Corollary \ref{eigenvalue.estimate.2}. From this we deduce that either $v_{(\bar{\xi},\bar{\varepsilon})} \geq 0$ almost everywhere, or 
\[\bigg ( \int_{\{v_{(\bar{\xi},\bar{\varepsilon})} < 0\}} |v_{(\bar{\xi},\bar{\varepsilon})}|^{\frac{2n}{n-2}} \bigg )^{\frac{2}{n}} \geq \frac{1}{2n(n-2) \, K}.\] 
On the other hand, we have 
\[\bigg ( \int_{\{v_{(\bar{\xi},\bar{\varepsilon})} < 0\}} |v_{(\bar{\xi},\bar{\varepsilon})}|^{\frac{2n}{n-2}} \bigg )^{\frac{n-2}{2n}} \leq \bigg ( \int_{\mathbb{R}^n} |v_{(\bar{\xi},\bar{\varepsilon})} - u_{(\bar{\xi},\bar{\varepsilon})}|^{\frac{2n}{n-2}} \bigg )^{\frac{n-2}{2n}} \leq C \, \alpha_1.\] 
Hence, if $\alpha_1$ is sufficiently small, then we have $v_{(\bar{\xi},\bar{\varepsilon})} \geq 0$ almost everywhere. \\

\section{An estimate for the energy of a "bubble"}

Throughout this paper, we fix a multi-linear form $W: \mathbb{R}^n \times \mathbb{R}^n \times \mathbb{R}^n \times \mathbb{R}^n \to \mathbb{R}$. We assume that $W_{ijkl}$ satisfy all the algebraic properties of the Weyl tensor. Moreover, we assume that some components of $W$ are non-zero, so that 
\[\sum_{i,j,k,l=1}^n (W_{ijkl} + W_{ilkj})^2 > 0.\] 
For abbreviation, we put 
\[H_{ik}(x) = \sum_{p,q=1}^n W_{ipkq} \, x_p \, x_q\] 
and 
\[\overline{H}_{ik}(x) = (1 - |x|^2) \, H_{ik}(x).\] 
It is easy to see that $H_{ik}(x)$ is trace-free, $\sum_{i=1}^n x_i \, H_{ik}(x) = 0$, and $\sum_{i=1}^n \partial_i H_{ik}(x) = 0$ for all $x \in \mathbb{R}^n$. \\

We consider a Riemannian metric of the form $g(x) = \exp(h(x))$, where $h(x)$ is a trace-free symmetric two-tensor on $\mathbb{R}^n$ satisfying $h(x) = 0$ for $|x| \geq 1$, 
\[|h(x)| + |\partial h(x)| + |\partial^2 h(x)| \leq \alpha_1\] 
for all $x \in \mathbb{R}^n$, and 
\[h_{ik}(x) = \mu \, (\lambda^2 - |x|^2) \, H_{ik}(x)\] 
for $|x| \leq \rho$. We assume that the parameters $\lambda$, $\mu$, and $\rho$ are chosen such that
$\mu \leq 1$ and $\lambda \leq \rho \leq 1$. Note that
$\sum_{i=1}^n x_i \, h_{ik}(x) = 0$ and $\sum_{i=1}^n \partial_i h_{ik}(x) = 0$ for $|x| \leq \rho$.

Given any pair $(\xi,\varepsilon) \in \mathbb{R}^n \times (0,\infty)$, there exists a unique function $v_{(\xi,\varepsilon)}$ such that $v_{(\xi,\varepsilon)} - u_{(\xi,\varepsilon)} \in \mathcal{E}_{(\xi,\varepsilon)}$ and 
\[\int_{\mathbb{R}^n} \Big ( \langle dv_{(\xi,\varepsilon)},d\psi \rangle_g + \frac{n-2}{4(n-1)} \, R_g \, v_{(\xi,\varepsilon)} \, \psi - n(n-2) \, |v_{(\xi,\varepsilon)}|^{\frac{4}{n-2}} \, v_{(\xi,\varepsilon)} \, \psi \Big ) = 0\] for all test functions $\psi \in \mathcal{E}_{(\xi,\varepsilon)}$ (see Proposition \ref{fixed.point.argument}). For abbreviation, let 
\[\Omega = \bigg \{ (\xi,\varepsilon) \in \mathbb{R}^n \times \mathbb{R}: |\xi| < 1, \, \frac{n-8}{3(n+4)} < \varepsilon^2 < \frac{2(n-8)}{3(n+4)} \bigg \}.\] 

\vspace{2mm}

\begin{proposition} 
\label{estimate.for.error.term}
For every pair $(\xi,\varepsilon) \in \lambda \, \Omega$, we have 
\begin{align*} 
&\Big \| \Delta_g u_{(\xi,\varepsilon)} 
- \frac{n-2}{4(n-1)} \, R_g \, u_{(\xi,\varepsilon)} 
+ n(n-2) \, u_{(\xi,\varepsilon)}^{\frac{n+2}{n-2}} \Big \|_{L^{\frac{2n}{n+2}}(\mathbb{R}^n)} \\ 
&\leq C \, \lambda^4 \, \mu + C \, \Big ( \frac{\lambda}{\rho} \Big )^{\frac{n-2}{2}} 
\end{align*}
and 
\begin{align*} 
&\Big \| \Delta_g u_{(\xi,\varepsilon)} 
- \frac{n-2}{4(n-1)} \, R_g \, u_{(\xi,\varepsilon)} 
+ n(n-2) \, u_{(\xi,\varepsilon)}^{\frac{n+2}{n-2}} \\ 
&\hspace{10mm} + \sum_{i,k=1}^n \mu \, (\lambda^2 - |x|^2) \, H_{ik}(x) \, \partial_i \partial_k u_{(\xi,\varepsilon)} \Big \|_{L^{\frac{2n}{n+2}}(\mathbb{R}^n)} \\ 
&\leq C \, \lambda^8 \, \mu^2 + C \, \Big ( \frac{\lambda}{\rho} \Big )^{\frac{n-2}{2}}. 
\end{align*}
\end{proposition} 

\textbf{Proof.} 
For abbreviation, we define two functions $A_1$ and $A_2$ by 
\[A_1 = \Delta_g u_{(\xi,\varepsilon)} 
- \frac{n-2}{4(n-1)} \, R_g \, u_{(\xi,\varepsilon)} 
+ n(n-2) \, u_{(\xi,\varepsilon)}^{\frac{n+2}{n-2}}\] 
and 
\[A_2 = \sum_{i,k=1}^n \mu \, (\lambda^2 - |x|^2) \, H_{ik}(x) \, \partial_i \partial_k u_{(\xi,\varepsilon)}.\] 
Using Proposition \ref{Taylor.expansion.of.scalar.curvature} and the 
identity $\sum_{i=1}^n \partial_i h_{ik}(x) = 0$, we obtain 
\[|R_g(x)| \leq C \, |h(x)|^2 \, |\partial^2 h(x)| + C \, |\partial h(x)|^2 \leq C \, \mu^2 \, (\lambda + |x|)^6\] 
for $|x| \leq \rho$. This implies 
\begin{align*} 
|A_1| &= \Big | \sum_{i,k=1}^n \partial_i \big [ (g^{ik} - \delta_{ik}) \, \partial_k u_{(\xi,\varepsilon)} \big ] - \frac{n-2}{4(n-1)} \, R_g \, u_{(\xi,\varepsilon)} \Big | \\ 
&\leq C \, \lambda^{\frac{n-2}{2}} \, \mu \, (\lambda + |x|)^{4-n}
\end{align*} 
and 
\begin{align*} 
|A_1 + A_2| &= \Big | \sum_{i,k=1}^n \partial_i \big [ (g^{ik} - \delta_{ik} + h_{ik}) \, \partial_k u_{(\xi,\varepsilon)} \big ] - \frac{n-2}{4(n-1)} \, R_g \, u_{(\xi,\varepsilon)} \Big | \\ 
&\leq C \, \lambda^{\frac{n-2}{2}} \, \mu^2 \, (\lambda + |x|)^{8-n}
\end{align*} 
for $|x| \leq \rho$. Hence, we obtain 
\[\|A_1\|_{L^{\frac{2n}{n+2}}(B_\rho(0))} \leq C \, \lambda^{\frac{n-2}{2}} \, \mu \, \bigg ( \int_{\mathbb{R}^n} (\lambda + |x|)^{-\frac{2n(n-4)}{n+2}} \bigg )^{\frac{n+2}{2n}} \leq C \, \lambda^4 \, \mu\] 
and 
\[\|A_1 + A_2\|_{L^{\frac{2n}{n+2}}(B_\rho(0))} \leq C \, \lambda^{\frac{n-2}{2}} \, \mu^2 \, \bigg ( \int_{\mathbb{R}^n} (\lambda + |x|)^{-\frac{2n(n-8)}{n+2}} \bigg )^{\frac{n+2}{2n}} \leq C \, \lambda^8 \, \mu^2.\] 
On the other hand, we have 
\[|A_1(x)| \leq C \, \lambda^{\frac{n-2}{2}} \, |x|^{-n}\] 
for $\rho \leq |x| \leq 1$ and 
\[|A_2(x)| \leq C \, \lambda^{\frac{n-2}{2}} \, \mu \, |x|^{4-n}\] 
for $|x| \geq \rho$. Since the function $A_1(x)$ vanishes for $|x| \geq 1$, 
we conclude that  
\[\|A_1\|_{L^{\frac{2n}{n+2}}(\mathbb{R}^n \setminus B_\rho(0))} \\ 
\leq C \, \lambda^{\frac{n-2}{2}} \, \bigg ( \int_{\mathbb{R}^n \setminus B_\rho(0)} |x|^{-\frac{2n^2}{n+2}} \bigg )^{\frac{n+2}{2n}} \leq C \, \Big ( \frac{\lambda}{\rho} \Big )^{\frac{n-2}{2}}\] 
and 
\[\|A_2\|_{L^{\frac{2n}{n+2}}(\mathbb{R}^n \setminus B_\rho(0))} \leq C \, \lambda^{\frac{n-2}{2}} \, \mu \, \bigg ( \int_{\mathbb{R}^n \setminus B_\rho(0)} |x|^{-\frac{2n(n-4)}{n+2}} \bigg )^{\frac{n+2}{2n}} \leq C \, \rho^4 \, \mu \, \Big ( \frac{\lambda}{\rho} \Big )^{\frac{n-2}{2}}.\] 
Putting these facts together, the assertion follows. \\

\begin{corollary} 
\label{estimate.for.v.1}
The function $v_{(\xi,\varepsilon)} - u_{(\xi,\varepsilon)}$ satisfies the estimate 
\[\|v_{(\xi,\varepsilon)} - u_{(\xi,\varepsilon)}\|_{L^{\frac{2n}{n-2}}(\mathbb{R}^n)} 
\leq C \, \lambda^4 \, \mu + C \, \Big ( \frac{\lambda}{\rho} \Big )^{\frac{n-2}{2}}\] 
for $(\xi,\varepsilon) \in \lambda \, \Omega$.
\end{corollary}

\textbf{Proof.} 
It follows from Proposition \ref{fixed.point.argument} that 
\begin{align*} 
&\|v_{(\xi,\varepsilon)} - u_{(\xi,\varepsilon)}\|_{L^{\frac{2n}{n-2}}(\mathbb{R}^n)} \\ 
&\leq C \, \Big \| \Delta_g u_{(\xi,\varepsilon)} 
- \frac{n-2}{4(n-1)} \, R_g \, u_{(\xi,\varepsilon)} 
+ n(n-2) \, u_{(\xi,\varepsilon)}^{\frac{n+2}{n-2}} \Big \|_{L^{\frac{2n}{n+2}}(\mathbb{R}^n)}, 
\end{align*} 
where $C$ is a constant that depends only on $n$. Hence, the assertion follows from Proposition \ref{estimate.for.error.term}. \\

We now prove a more refined estimate for the difference $v_{(\xi,\varepsilon)} - u_{(\xi,\varepsilon)}$. 
Using Proposition \ref{linearized.operator} with $h = 0$, we conclude that there exists a unique function $w_{(\xi,\varepsilon)} \in \mathcal{E}_{(\xi,\varepsilon)}$ such that 
\begin{align} 
\label{def.w} 
&\int_{\mathbb{R}^n} \Big ( \langle dw_{(\xi,\varepsilon)},d\psi \rangle - n(n+2) \, u_{(\xi,\varepsilon)}^{\frac{4}{n-2}} \, w_{(\xi,\varepsilon)} \, \psi \Big ) \notag \\ 
&= -\int_{\mathbb{R}^n} \sum_{i,k=1}^n \mu \, (\lambda^2 - |x|^2) \, H_{ik}(x) \, \partial_i \partial_k u_{(\xi,\varepsilon)} \, \psi 
\end{align} 
for all test functions $\psi \in \mathcal{E}_{(\xi,\varepsilon)}$. 

\begin{proposition} 
\label{properties.of.w}
The function $w_{(\xi,\varepsilon)}$ is smooth. Moreover, if $(\xi,\varepsilon) \in \lambda \, \Omega$, then we have 
\begin{align*} 
&|w_{(\xi,\varepsilon)}(x)| \leq C \, \lambda^{\frac{n-2}{2}} \, \mu \, (\lambda + |x|)^{6-n} \\ 
&|\partial w_{(\xi,\varepsilon)}(x)| \leq C \, \lambda^{\frac{n-2}{2}} \, \mu \, (\lambda + |x|)^{5-n} \\ 
&|\partial^2 w_{(\xi,\varepsilon)}(x)| \leq C \, \lambda^{\frac{n-2}{2}} \, \mu \, (\lambda + |x|)^{4-n} 
\end{align*}
for all $x \in \mathbb{R}^n$.
\end{proposition}

\textbf{Proof.} 
Let $\varphi_{(\xi,\varepsilon,k)}$ be the functions defined in Section 2. We can find real numbers $b_k(\xi,\varepsilon)$, $k = 0,1,\hdots,n$, such that 
\begin{align*} 
&\int_{\mathbb{R}^n} \Big ( \langle dw_{(\xi,\varepsilon)},d\psi \rangle - n(n+2) \, u_{(\xi,\varepsilon)}^{\frac{4}{n-2}} \, w_{(\xi,\varepsilon)} \, \psi \Big ) \notag \\ 
&= -\int_{\mathbb{R}^n} \sum_{i,k=1}^n \mu \, (\lambda^2 - |x|^2) \, H_{ik}(x) \, \partial_i \partial_k u_{(\xi,\varepsilon)} \, \psi + \sum_{k=0}^n b_k(\xi,\varepsilon) \, \int_{\mathbb{R}^n} \varphi_{(\xi,\varepsilon,k)} \, \psi 
\end{align*} 
for all test functions $\psi \in \mathcal{E}$. It follows from standard elliptic regularity theory that $w_{(\xi,\varepsilon)}$ is smooth. 

In the next step, we establish quantitative estimates for $w_{(\xi,\varepsilon)}$. To that end, we consider a pair $(\xi,\varepsilon) \in \lambda \, \Omega$. A straightforward calculation yields 
\begin{equation} 
\label{estimate}
\Big \| \sum_{i,k=1}^n \mu \, (\lambda^2 - |x|^2) \, H_{ik}(x) \, \partial_i \partial_k u_{(\xi,\varepsilon)} \Big \|_{L^{\frac{2n}{n+2}}(\mathbb{R}^n)} \leq C \, \lambda^4 \, \mu. 
\end{equation}
From this we deduce that $\|w_{(\xi,\varepsilon)}\|_{L^{\frac{2n}{n-2}}(\mathbb{R}^n)} \leq C \, \lambda^4 \, \mu$ and $\sum_{k=0}^n |b_k(\xi,\varepsilon)| \leq C \, \lambda^4 \, \mu$. 
This implies 
\begin{align*} 
&\big | \Delta w_{(\xi,\varepsilon)} + n(n+2) \, u_{(\xi,\varepsilon)}^{\frac{4}{n-2}} \, w_{(\xi,\varepsilon)} \big | \\ 
&= \bigg | \sum_{i,k=1}^n \mu \, (\lambda^2 - |x|^2) \, H_{ik}(x) \, \partial_i \partial_k u_{(\xi,\varepsilon)} - \sum_{k=0}^n b_k(\xi,\varepsilon) \, \varphi_{(\xi,\varepsilon,k)} \bigg | \\ 
&\leq C \, \lambda^{\frac{n-2}{2}} \, \mu \, (\lambda + |x|)^{4-n}
\end{align*} 
for all $x \in \mathbb{R}^n$. We claim that 
\[\sup_{x \in \mathbb{R}^n} (\lambda + |x|)^{\frac{n-2}{2}} \, |w_{(\xi,\varepsilon)}(x)| \leq C \, \lambda^4 \, \mu.\] 
To show this, we fix a point $x_0 \in \mathbb{R}^n$ and put $r = \frac{1}{2} \, (\lambda + |x_0|)$. Clearly, 
$\lambda + |x| \geq r$ for all $x \in B_r(x_0)$. This implies \[u_{(\xi,\varepsilon)}(x)^{\frac{4}{n-2}} \leq C \, r^{-2}\] and 
\[\big | \Delta w_{(\xi,\varepsilon)} + n(n+2) \, u_{(\xi,\varepsilon)}^{\frac{4}{n-2}} \, w_{(\xi,\varepsilon)} \big | \leq C \, \lambda^{\frac{n-2}{2}} \, \mu \, r^{4-n}\] 
for all $x \in B_r(x_0)$. Using standard interior estimates, we obtain 
\begin{align*} 
r^{\frac{n-2}{2}} \, |w_{(\xi,\varepsilon)}(x_0)| &\leq C \, \|w_{(\xi,\varepsilon)}\|_{L^{\frac{2n}{n-2}}(B_r(x_0))} \\ 
&+ C \, r^{\frac{n+2}{2}} \, \big \| \Delta w_{(\xi,\varepsilon)} + n(n+2) \, u_{(\xi,\varepsilon)}^{\frac{4}{n-2}} \, w_{(\xi,\varepsilon)} \big \|_{L^\infty(B_r(x_0))} \\ 
&\leq C \, \lambda^4 \, \mu + C \, \lambda^{\frac{n-2}{2}} \, \mu \, r^{-\frac{n-10}{2}} \\ 
&\leq C \, \lambda^4 \, \mu. 
\end{align*} 
Thus, we conclude that 
\[\sup_{x \in \mathbb{R}^n} (\lambda + |x|)^{\frac{n-2}{2}} \, |w_{(\xi,\varepsilon)}(x)| \leq C \, \lambda^4 \, \mu,\] as claimed. Since $\sup_{x \in \mathbb{R}^n} |x|^{\frac{n-2}{2}} \, |w_{(\xi,\varepsilon)}(x)| < \infty$, we can express the function $w_{(\xi,\varepsilon)}$ in the form 
\begin{equation} 
\label{convolution.formula}
w_{(\xi,\varepsilon)}(x) = -\frac{1}{(n-2) \, |S^{n-1}|} \int_{\mathbb{R}^n} |x - y|^{2-n} \, \Delta w_{(\xi,\varepsilon)}(y) \, dy 
\end{equation}
for all $x \in \mathbb{R}^n$. 

We can now use a bootstrap argument to prove the desired estimate for $w_{(\xi,\varepsilon)}$. It follows from (\ref{convolution.formula}) that 
\[\sup_{x \in \mathbb{R}^n} (\lambda + |x|)^\beta \, |w_{(\xi,\varepsilon)}(x)| \leq C \, \sup_{x \in \mathbb{R}^n} (\lambda + |x|)^{\beta+2} \, |\Delta w_{(\xi,\varepsilon)}(x)|\] 
for all $0 < \beta < n-2$. Since 
\begin{align*} 
|\Delta w_{(\xi,\varepsilon)}(x)| &\leq n(n+2) \, u_{(\xi,\varepsilon)}(x)^{\frac{4}{n-2}} \, |w_{(\xi,\varepsilon)}(x)| \\ 
&+ C \, \lambda^{\frac{n-2}{2}} \, \mu \, (\lambda + |x|)^{4-n} 
\end{align*} 
for all $x \in \mathbb{R}^n$, we conclude that 
\begin{align*} 
\sup_{x \in \mathbb{R}^n} (\lambda + |x|)^\beta \, |w_{(\xi,\varepsilon)}(x)| &\leq C \, \lambda^2 \, \sup_{x \in \mathbb{R}^n} (\lambda + |x|)^{\beta-2} \, |w_{(\xi,\varepsilon)}(x)| \\ &+ C \, \lambda^{\beta-\frac{n-10}{2}} \, \mu 
\end{align*} 
for all $0 < \beta \leq n-6$. Iterating this inequality, we obtain 
\[\sup_{x \in \mathbb{R}^n} (\lambda + |x|)^{n-6} \, |w_{(\xi,\varepsilon)}(x)| \leq C \, \lambda^{\frac{n-2}{2}} \, \mu.\] 
The estimates for the first and second derivatives of $w_{(\xi,\varepsilon)}$ follow now from standard interior estimates. \\

\begin{corollary} 
\label{estimate.for.v.2}
The function $v_{(\xi,\varepsilon)} - u_{(\xi,\varepsilon)} - w_{(\xi,\varepsilon)}$ satisfies the estimate 
\[\|v_{(\xi,\varepsilon)} - u_{(\xi,\varepsilon)} - w_{(\xi,\varepsilon)}\|_{L^{\frac{2n}{n-2}}(\mathbb{R}^n)} 
\leq C \, \lambda^{\frac{4(n+2)}{n-2}} \, \mu^{\frac{n+2}{n-2}} + C \, \Big ( \frac{\lambda}{\rho} \Big )^{\frac{n-2}{2}}\] 
for $(\xi,\varepsilon) \in \lambda \, \Omega$.
\end{corollary}

\textbf{Proof.} 
Consider the functions 
\[B_1 = \sum_{i,k=1}^n \partial_i \big [ (g^{ik} - \delta_{ik}) \, \partial_k w_{(\xi,\varepsilon)} \big ] - \frac{n-2}{4(n-1)} \, R_g \, w_{(\xi,\varepsilon)}\] 
and 
\[B_2 = \sum_{i,k=1}^n \mu \, (\lambda^2 - |x|^2) \, H_{ik}(x) \, \partial_i \partial_k u_{(\xi,\varepsilon)}.\] 
Using (\ref{def.w}), we obtain 
\begin{align*} 
&\int_{\mathbb{R}^n} \Big ( \langle dw_{(\xi,\varepsilon)},d\psi \rangle_g 
+ \frac{n-2}{4(n-1)} \, R_g \, w_{(\xi,\varepsilon)} \, \psi - n(n+2) \, u_{(\xi,\varepsilon)}^{\frac{4}{n-2}} \, w_{(\xi,\varepsilon)} \, \psi \Big ) \\ 
&= -\int_{\mathbb{R}^n} (B_1 + B_2) \, \psi 
\end{align*} 
for all functions $\psi \in \mathcal{E}_{(\xi,\varepsilon)}$. Since $w_{(\xi,\varepsilon)} \in \mathcal{E}_{(\xi,\varepsilon)}$, it follows that 
\[w_{(\xi,\varepsilon)} = -G_{(\xi,\varepsilon)}(B_1 + B_2).\] 
Moreover, we have 
\[v_{(\xi,\varepsilon)} - u_{(\xi,\varepsilon)} = G_{(\xi,\varepsilon)} \big ( B_3 + n(n-2) \, B_4 \big ),\] 
where 
\[B_3 = \Delta_g u_{(\xi,\varepsilon)} - \frac{n-2}{4(n-1)} \, R_g \, u_{(\xi,\varepsilon)} + n(n-2) \, u_{(\xi,\varepsilon)}^{\frac{n+2}{n-2}}\] 
and 
\[B_4 = |v_{(\xi,\varepsilon)}|^{\frac{4}{n-2}} \, v_{(\xi,\varepsilon)} - u_{(\xi,\varepsilon)}^{\frac{n+2}{n-2}}  - \frac{n+2}{n-2} \, u_{(\xi,\varepsilon)}^{\frac{4}{n-2}} \, (v_{(\xi,\varepsilon)} - u_{(\xi,\varepsilon)}).\] 
Thus, we conclude that 
\[v_{(\xi,\varepsilon)} - u_{(\xi,\varepsilon)} - w_{(\xi,\varepsilon)} = G_{(\xi,\varepsilon)} \big ( B_1 + B_2 + B_3 + n(n-2) \, B_4 \big ),\] 
where $G_{(\xi,\varepsilon)}: L^{\frac{2n}{n+2}}(\mathbb{R}^n) \to \mathcal{E}_{(\xi,\varepsilon)}$ denotes the solution operator constructed in Proposition \ref{linearized.operator}. In particular, we have 
\[\|v_{(\xi,\varepsilon)} - u_{(\xi,\varepsilon)} - w_{(\xi,\varepsilon)}\|_{L^{\frac{2n}{n-2}}(\mathbb{R}^n)} 
\leq C \, \big \| B_1 + B_2 + B_3 + n(n-2) \, B_4 \big \|_{L^{\frac{2n}{n+2}}(\mathbb{R}^n)}\] 
by Proposition \ref{linearized.operator}. Using Proposition \ref{properties.of.w}, we obtain 
\[|B_1(x)| \leq C \, \lambda^{\frac{n-2}{2}} \, \mu^2 \, (\lambda + |x|)^{8-n}\] 
for $|x| \leq \rho$ and 
\[|B_1(x)| \leq C \, \lambda^{\frac{n-2}{2}} \, \mu \, |x|^{4-n}\] 
for $\rho \leq |x| \leq 1$. Since the function $B_1(x)$ vanishes for $|x| \geq 1$, 
we conclude that 
\[\|B_1\|_{L^{\frac{2n}{n+2}}(\mathbb{R}^n)} \leq C \, \lambda^8 \, \mu^2 + C \, \rho^4 \, \mu \, \Big ( \frac{\lambda}{\rho} \Big )^{\frac{n-2}{2}}.\] 
Moreover, we have 
\[\|B_2 + B_3\|_{L^{\frac{2n}{n+2}}(\mathbb{R}^n)} \leq C \, \lambda^8 \, \mu^2 + C \, \Big ( \frac{\lambda}{\rho} \Big )^{\frac{n-2}{2}}\] 
by Proposition \ref{estimate.for.error.term}. Finally, the function $B_4$ satisfies a pointwise estimate of the form 
\[|B_4| \leq C \, |v_{(\xi,\varepsilon)} - u_{(\xi,\varepsilon)}|^{\frac{n+2}{n-2}},\] 
where $C$ is a constant that depends only on $n$. Hence, it follows from Corollary \ref{estimate.for.v.1} that 
\begin{align*} 
\|B_4\|_{L^{\frac{2n}{n+2}}(\mathbb{R}^n)} 
&\leq C \, \|v_{(\xi,\varepsilon)} - u_{(\xi,\varepsilon)}\|_{L^{\frac{2n}{n-2}}(\mathbb{R}^n)}^{\frac{n+2}{n-2}} \\ 
&\leq C \, \lambda^{\frac{4(n+2)}{n-2}} \, \mu^{\frac{n+2}{n-2}} + C \, \Big ( \frac{\lambda}{\rho} \Big )^{\frac{n+2}{2}}. 
\end{align*} 
Putting these facts together, we obtain 
\[\|v_{(\xi,\varepsilon)} - u_{(\xi,\varepsilon)} - w_{(\xi,\varepsilon)}\|_{L^{\frac{2n}{n-2}}(\mathbb{R}^n)} \leq C \, \lambda^{\frac{4(n+2)}{n-2}} \, \mu^{\frac{n+2}{n-2}} + C \, \Big ( \frac{\lambda}{\rho} \Big )^{\frac{n-2}{2}},\] 
as claimed. \\

\begin{proposition} 
\label{term.1}
We have 
\begin{align*} 
&\bigg | \int_{\mathbb{R}^n} \Big ( |dv_{(\xi,\varepsilon)}|_g^2 
- |du_{(\xi,\varepsilon)}|_g^2 + \frac{n-2}{4(n-1)} \, R_g \, 
(v_{(\xi,\varepsilon)}^2 - u_{(\xi,\varepsilon)}^2) \Big ) \\ 
&\hspace{10mm} + \int_{\mathbb{R}^n} n(n-2) \, (|v_{(\xi,\varepsilon)}|^{\frac{4}{n-2}} - u_{(\xi,\varepsilon)}^{\frac{4}{n-2}}) \, u_{(\xi,\varepsilon)} \, v_{(\xi,\varepsilon)} \\ 
&\hspace{10mm} - \int_{\mathbb{R}^n} n(n-2) \, (|v_{(\xi,\varepsilon)}|^{\frac{2n}{n-2}} - u_{(\xi,\varepsilon)}^{\frac{2n}{n-2}}) \\ 
&\hspace{10mm} - \int_{\mathbb{R}^n} \sum_{i,k=1}^n \mu \, (\lambda^2 - |x|^2) \, H_{ik}(x) \, \partial_i \partial_k u_{(\xi,\varepsilon)} \, w_{(\xi,\varepsilon)} \bigg | \\ 
&\leq C \, \lambda^{\frac{8n}{n-2}} \, \mu^{\frac{2n}{n-2}} + C \, \lambda^4 \, \mu \, \Big ( \frac{\lambda}{\rho} \Big )^{\frac{n-2}{2}} + C \, \Big ( \frac{\lambda}{\rho} \Big )^{n-2} 
\end{align*}
for $(\xi,\varepsilon) \in \lambda \, \Omega$.
\end{proposition}

\textbf{Proof.} 
Using Proposition \ref{fixed.point.argument} with $\psi = v_{(\xi,\varepsilon)} - u_{(\xi,\varepsilon)}$, we obtain 
\begin{align*} 
&\int_{\mathbb{R}^n} \Big ( |dv_{(\xi,\varepsilon)}|_g^2 - \langle du_{(\xi,\varepsilon)},dv_{(\xi,\varepsilon)} \rangle_g 
+ \frac{n-2}{4(n-1)} \, R_g \, v_{(\xi,\varepsilon)} \, (v_{(\xi,\varepsilon)} - u_{(\xi,\varepsilon)}) \Big ) \\ 
&\hspace{10mm} - \int_{\mathbb{R}^n} n(n-2) \, |v_{(\xi,\varepsilon)}|^{\frac{4}{n-2}} \, v_{(\xi,\varepsilon)} \, (v_{(\xi,\varepsilon)} - u_{(\xi,\varepsilon)}) = 0. 
\end{align*} 
Moreover, it follows from Proposition \ref{estimate.for.error.term} and Corollary \ref{estimate.for.v.1} that 
\begin{align*} 
&\bigg | \int_{\mathbb{R}^n} \Big ( \langle du_{(\xi,\varepsilon)},dv_{(\xi,\varepsilon)} \rangle_g 
- |du_{(\xi,\varepsilon)}|_g^2 + \frac{n-2}{4(n-1)} \, R_g \, u_{(\xi,\varepsilon)} \, (v_{(\xi,\varepsilon)} - u_{(\xi,\varepsilon)}) \Big ) \\ 
&\hspace{10mm} - \int_{\mathbb{R}^n} n(n-2) \, u_{(\xi,\varepsilon)}^{\frac{n+2}{n-2}} \, (v_{(\xi,\varepsilon)} - u_{(\xi,\varepsilon)}) \\ 
&\hspace{10mm} - \int_{\mathbb{R}^n} \sum_{i,k=1}^n \mu \, (\lambda^2 - |x|^2) \, H_{ik}(x) \, \partial_i \partial_k u_{(\xi,\varepsilon)} \, (v_{(\xi,\varepsilon)} - u_{(\xi,\varepsilon)}) \bigg | \\ 
&\leq \Big \| \Delta_g u_{(\xi,\varepsilon)} 
- \frac{n-2}{4(n-1)} \, R_g \, u_{(\xi,\varepsilon)} 
+ n(n-2) \, u_{(\xi,\varepsilon)}^{\frac{n+2}{n-2}} \\ 
&\hspace{10mm} + \sum_{i,k=1}^n \mu \, (\lambda^2 - |x|^2) \, H_{ik}(x) \, \partial_i \partial_k u_{(\xi,\varepsilon)} \Big \|_{L^{\frac{2n}{n+2}}(\mathbb{R}^n)} \\ 
&\hspace{5mm} \cdot \|v_{(\xi,\varepsilon)} - u_{(\xi,\varepsilon)}\|_{L^{\frac{2n}{n-2}}(\mathbb{R}^n)} \\ 
&\leq C \, \lambda^{12} \, \mu^3 + C \, \lambda^4 \, \mu \, \Big ( \frac{\lambda}{\rho} \Big )^{\frac{n-2}{2}} + C \, \Big ( \frac{\lambda}{\rho} \Big )^{n-2}. 
\end{align*}
Finally, we have 
\begin{align*} 
&\bigg | \int_{\mathbb{R}^n} \sum_{i,k=1}^n \mu \, (\lambda^2 - |x|^2) \, H_{ik}(x) \, \partial_i \partial_k u_{(\xi,\varepsilon)} \, (v_{(\xi,\varepsilon)} - u_{(\xi,\varepsilon)} - w_{(\xi,\varepsilon)}) \bigg | \\ 
&\leq C \, \lambda^4 \, \mu \, \|v_{(\xi,\varepsilon)} - u_{(\xi,\varepsilon)} - w_{(\xi,\varepsilon)}\|_{L^{\frac{2n}{n-2}}(\mathbb{R}^n)} \\ 
&\leq C \, \lambda^{\frac{8n}{n-2}} \, \mu^{\frac{2n}{n-2}} + C \, \lambda^4 \, \mu \, \Big ( \frac{\lambda}{\rho} \Big )^{\frac{n-2}{2}} 
\end{align*} 
by (\ref{estimate}) and Corollary \ref{estimate.for.v.2}. Putting these facts together, the assertion follows. \\

\begin{proposition} 
\label{term.2}
We have 
\begin{align*} 
&\bigg | \int_{\mathbb{R}^n} (|v_{(\xi,\varepsilon)}|^{\frac{4}{n-2}} - u_{(\xi,\varepsilon)}^{\frac{4}{n-2}}) \, u_{(\xi,\varepsilon)} \, v_{(\xi,\varepsilon)} - \frac{2}{n} \int_{\mathbb{R}^n} (|v_{(\xi,\varepsilon)}|^{\frac{2n}{n-2}} - u_{(\xi,\varepsilon)}^{\frac{2n}{n-2}}) \bigg | \\ 
&\leq C \, \lambda^{\frac{8n}{n-2}} \, \mu^{\frac{2n}{n-2}} + C \, \Big ( \frac{\lambda}{\rho} \Big )^n
\end{align*}
for $(\xi,\varepsilon) \in \lambda \, \Omega$.
\end{proposition} 

\textbf{Proof.} 
We have the pointwise estimate 
\begin{align*} 
&\Big | (|v_{(\xi,\varepsilon)}|^{\frac{4}{n-2}} - u_{(\xi,\varepsilon)}^{\frac{4}{n-2}}) \, u_{(\xi,\varepsilon)} \, v_{(\xi,\varepsilon)} - \frac{2}{n} \, (|v_{(\xi,\varepsilon)}|^{\frac{2n}{n-2}} - u_{(\xi,\varepsilon)}^{\frac{2n}{n-2}}) \Big | \\ 
&\leq C \, |v_{(\xi,\varepsilon)} - u_{(\xi,\varepsilon)}|^{\frac{2n}{n-2}}, 
\end{align*} 
where $C$ is a constant that depends only on $n$. This implies
\begin{align*} 
&\bigg | \int_{\mathbb{R}^n} (|v_{(\xi,\varepsilon)}|^{\frac{4}{n-2}} - u_{(\xi,\varepsilon)}^{\frac{4}{n-2}}) \, u_{(\xi,\varepsilon)} \, v_{(\xi,\varepsilon)} - \frac{2}{n} \int_{\mathbb{R}^n} (|v_{(\xi,\varepsilon)}|^{\frac{2n}{n-2}} - u_{(\xi,\varepsilon)}^{\frac{2n}{n-2}}) \bigg | \\ 
&\leq C \, \|v_{(\xi,\varepsilon)} - u_{(\xi,\varepsilon)}\|_{L^{\frac{2n}{n-2}}(\mathbb{R}^n)}^{\frac{2n}{n-2}} \\ 
&\leq C \, \lambda^{\frac{8n}{n-2}} \, \mu^{\frac{2n}{n-2}} + C \, \Big ( \frac{\lambda}{\rho} \Big )^n 
\end{align*} 
by Corollary \ref{estimate.for.v.1}. \\

\begin{proposition} 
\label{term.3}
We have 
\begin{align*} 
&\bigg | \int_{\mathbb{R}^n} \Big ( |du_{(\xi,\varepsilon)}|_g^2 + \frac{n-2}{4(n-1)} \, R_g \, u_{(\xi,\varepsilon)}^2 - n(n-2) \, u_{(\xi,\varepsilon)}^{\frac{2n}{n-2}} \Big ) \\ 
&\hspace{10mm} - \int_{B_\rho(0)} \frac{1}{2} \, \sum_{i,k,l=1}^n h_{il} \, h_{kl} \, \partial_i u_{(\xi,\varepsilon)} \, \partial_k u_{(\xi,\varepsilon)} \\ 
&\hspace{10mm} + \int_{B_\rho(0)} \frac{n-2}{16(n-1)} \, \sum_{i,k,l=1}^n (\partial_l h_{ik})^2 \, 
u_{(\xi,\varepsilon)}^2 \bigg | \\ 
&\leq C \, \lambda^{12} \, \mu^3 + C \, \Big ( \frac{\lambda}{\rho} \Big )^{n-2} 
\end{align*} 
for all $(\xi,\varepsilon) \in \lambda \, \Omega$.
\end{proposition}

\textbf{Proof.} 
Note that 
\begin{align*} 
&\Big | g^{ik}(x) - \delta_{ik} + h_{ik}(x) - \frac{1}{2} \, \sum_{l=1}^n h_{il}(x) \, h_{kl}(x) \Big | \\ 
&\leq C \, |h(x)|^3 \leq C \, \mu^3 \, (\lambda + |x|)^{12} 
\end{align*} 
for $|x| \leq \rho$. This implies 
\begin{align*} 
&\bigg | \int_{\mathbb{R}^n} \big ( |du_{(\xi,\varepsilon)}|_g^2 - |du_{(\xi,\varepsilon)}|^2 \big ) + \int_{\mathbb{R}^n} \sum_{i,k=1}^n h_{ik} \, \partial_i u_{(\xi,\varepsilon)} \, 
\partial_k u_{(\xi,\varepsilon)} \\ 
&\hspace{10mm} - \int_{B_\rho(0)} \frac{1}{2} \sum_{i,k,l=1}^n h_{il} \, h_{kl} \, \partial_i u_{(\xi,\varepsilon)} \, \partial_k u_{(\xi,\varepsilon)} \bigg | \\ 
&\leq C \, \lambda^{n-2} \, \mu^3 \, \int_{B_\rho(0)} (\lambda + |x|)^{14-2n} 
+ C \, \lambda^{n-2} \, \int_{\mathbb{R}^n \setminus B_\rho(0)} (\lambda + |x|)^{2-2n} \\ 
&\leq C \, \lambda^{12} \, \mu^3 + C \, \Big ( \frac{\lambda}{\rho} \Big )^{n-2}. 
\end{align*} 
By Proposition \ref{Taylor.expansion.of.scalar.curvature}, the scalar curvature of $g$ satisfies the estimate 
\begin{align*} 
&\Big | R_g(x) + \frac{1}{4} \sum_{i,k,l=1}^n (\partial_l h_{ik}(x))^2 \Big | \\ 
&\leq C \, |h(x)|^2 \, |\partial^2 h(x)| + C \, |h(x)| \, |\partial h(x)|^2 \\ 
&\leq C \, \mu^3 \, (\lambda + |x|)^{10} 
\end{align*} 
for $|x| \leq \rho$. This implies 
\begin{align*} 
&\bigg | \int_{\mathbb{R}^n} R_g \, u_{(\xi,\varepsilon)}^2 + \int_{B_\rho(0)} \frac{1}{4} \sum_{i,k,l=1}^n (\partial_l h_{ik})^2 \, u_{(\xi,\varepsilon)}^2 \bigg | \\ 
&\leq C \, \lambda^{12} \, \mu^3 \, \int_{B_\rho(0)} (\lambda + |x|)^{14-2n} 
+ C \, \lambda^{n-2} \, \int_{\mathbb{R}^n \setminus B_\rho(0)} (\lambda + |x|)^{4-2n} \\ 
&\leq C \, \lambda^{12} \, \mu^3 + C \, \rho^2 \, \Big ( \frac{\lambda}{\rho} \Big )^{n-2}. 
\end{align*}
At this point, we use the formula 
\begin{align*} 
&\partial_i u_{(\xi,\varepsilon)} \, \partial_k u_{(\xi,\varepsilon)} - \frac{n-2}{4(n-1)} \, \partial_i \partial_k (u_{(\xi,\varepsilon)}^2) \\ 
&= \frac{1}{n} \, \Big ( |du_{(\xi,\varepsilon)}|^2 - \frac{n-2}{4(n-1)} \, \Delta (u_{(\xi,\varepsilon)}^2) \Big ) \, \delta_{ik}. 
\end{align*}
Since $h_{ik}$ is trace-free, we obtain 
\[\sum_{i,k=1}^n h_{ik} \, \partial_i u_{(\xi,\varepsilon)} \, \partial_k u_{(\xi,\varepsilon)} = \frac{n-2}{4(n-1)} \sum_{i,k=1}^n h_{ik} \, \partial_i \partial_k (u_{(\xi,\varepsilon)}^2),\] 
hence 
\[\int_{\mathbb{R}^n} \sum_{i,k=1}^n h_{ik} \, \partial_i u_{(\xi,\varepsilon)} \, \partial_k u_{(\xi,\varepsilon)} = \int_{\mathbb{R}^n} \frac{n-2}{4(n-1)} \sum_{i,k=1}^n \partial_i \partial_k h_{ik} \, u_{(\xi,\varepsilon)}^2.\]
Since $\sum_{i=1}^n \partial_i h_{ik}(x) = 0$ for $|x| \leq \rho$, it follows that 
\[\bigg | \int_{\mathbb{R}^n} \sum_{i,k=1}^n h_{ik} \, \partial_i u_{(\xi,\varepsilon)} \, \partial_k u_{(\xi,\varepsilon)} \bigg | 
\leq C \int_{\mathbb{R}^n \setminus B_\rho(0)} u_{(\xi,\varepsilon)}^2 \leq C \, \rho^2 \, \Big ( \frac{\lambda}{\rho} \Big )^{n-2}.\] 
Putting these facts together, the assertion follows. \\

\begin{corollary} 
\label{key.estimate}
The function $\mathcal{F}_g(\xi,\varepsilon)$ satisfies the estimate 
\begin{align*} 
&\bigg | \mathcal{F}_g(\xi,\varepsilon) - \int_{B_\rho(0)} \frac{1}{2} \, \sum_{i,k,l=1}^n h_{il} \, h_{kl} \, \partial_i u_{(\xi,\varepsilon)} \, \partial_k u_{(\xi,\varepsilon)} \\ 
&\hspace{10mm} + \int_{B_\rho(0)} \frac{n-2}{16(n-1)} \, \sum_{i,k,l=1}^n (\partial_l h_{ik})^2 \, 
u_{(\xi,\varepsilon)}^2 \\ 
&\hspace{10mm} - \int_{\mathbb{R}^n} \sum_{i,k=1}^n 
\mu \, (\lambda^2 - |x|^2) \, H_{ik}(x) \, \partial_i\partial_k u_{(\xi,\varepsilon)} \, w_{(\xi,\varepsilon)} \bigg | \\ 
&\leq C \, \lambda^{\frac{8n}{n-2}} \, \mu^{\frac{2n}{n-2}} 
+ C \, \lambda^4 \, \mu \, \Big ( \frac{\lambda}{\rho} \Big )^{\frac{n-2}{2}} + C \, \Big ( \frac{\lambda}{\rho} \Big )^{n-2} 
\end{align*} 
for $(\xi,\varepsilon) \in \lambda \, \Omega$.
\end{corollary}

\textbf{Proof.} This follows by combining Proposition \ref{term.1}, Proposition \ref{term.2}, and Proposition \ref{term.3}. \\

\section{Finding a critical point of an auxiliary function}

We define a function $F: \mathbb{R}^n \times (0,\infty) \to \mathbb{R}$ as follows: given any pair $(\xi,\varepsilon) \in \mathbb{R}^n \times (0,\infty)$, we define 
\begin{align*} 
F(\xi,\varepsilon) 
&= \int_{\mathbb{R}^n} \frac{1}{2} \sum_{i,k,l=1}^n \overline{H}_{il}(x) \, \overline{H}_{kl}(x) \, \partial_i u_{(\xi,\varepsilon)}(x) \, \partial_k u_{(\xi,\varepsilon)}(x) \\ 
&- \int_{\mathbb{R}^n} \frac{n-2}{16(n-1)} \, \sum_{i,k,l=1}^n (\partial_l \overline{H}_{ik}(x))^2 \, u_{(\xi,\varepsilon)}(x)^2 \\ 
&+ \int_{\mathbb{R}^n} \sum_{i,k=1}^n \overline{H}_{ik}(x) \, \partial_i \partial_k u_{(\xi,\varepsilon)}(x) \, z_{(\xi,\varepsilon)}(x),
\end{align*} 
where $z_{(\xi,\varepsilon)} \in \mathcal{E}_{(\xi,\varepsilon)}$ satisfies the relation 
\begin{align*} 
&\int_{\mathbb{R}^n} \Big ( \langle dz_{(\xi,\varepsilon)},d\psi \rangle - n(n+2) \, u_{(\xi,\varepsilon)}(x)^{\frac{4}{n-2}} \, z_{(\xi,\varepsilon)} \, \psi \Big ) \\ 
&= -\int_{\mathbb{R}^n} \sum_{i,k=1}^n \overline{H}_{ik} \, \partial_i \partial_k u_{(\xi,\varepsilon)} \, \psi 
\end{align*} 
for all test functions $\psi \in \mathcal{E}_{(\xi,\varepsilon)}$. Our goal in this section is to show that the function $F(\xi,\varepsilon)$ has a critical point.

\begin{proposition} 
\label{symmetry} 
The function $F(\xi,\varepsilon)$ satisfies $F(\xi,\varepsilon) = F(-\xi,\varepsilon)$ for all $(\xi,\varepsilon) \in \mathbb{R}^n \times (0,\infty)$. Consequently, we have $\frac{\partial}{\partial \xi_p} F(0,\varepsilon) = 0$ and $\frac{\partial^2}{\partial \varepsilon \, \partial \xi_p} F(0,\varepsilon) = 0$ for all $\varepsilon > 0$ and $p = 1, \hdots, n$.
\end{proposition}

\textbf{Proof.} 
This follows immediately from the relation $\overline{H}_{ik}(-x) = \overline{H}_{ik}(x)$. \\

\begin{proposition} 
\label{integral.identity.1}
We have 
\begin{align*} 
&\int_{\partial B_r(0)} \sum_{i,k,l=1}^n (\partial_l H_{ik}(x))^2 \, x_p \, x_q \\ 
&= \frac{2}{n(n+2)} \, |S^{n-1}| \, \sum_{i,k,l=1}^n (W_{ipkl} + W_{ilkp}) \, (W_{iqkl} + W_{ilkq}) \, r^{n+3} \\ 
&+ \frac{1}{n(n+2)} \, |S^{n-1}| \, \sum_{i,j,k,l=1}^n (W_{ijkl} + W_{ilkj})^2 \, \delta_{pq} \, r^{n+3} 
\end{align*} 
and 
\begin{align*} 
&\int_{\partial B_r(0)} \sum_{i,k=1}^n H_{ik}(x)^2 \, x_p \, x_q \\ 
&= \frac{2}{n(n+2)(n+4)} \, |S^{n-1}| \, \sum_{i,k,l=1}^n (W_{ipkl} + W_{ilkp}) \, (W_{iqkl} + W_{ilkq}) \, r^{n+5} \\ 
&+ \frac{1}{2n(n+2)(n+4)} \, |S^{n-1}| \, \sum_{i,j,k,l=1}^n (W_{ijkl} + W_{ilkj})^2 \, \delta_{pq} \, r^{n+5}. 
\end{align*} 
\end{proposition}

\textbf{Proof.} 
By definition of $H_{ik}(x)$, we have 
\begin{align*} 
&\int_{\partial B_r(0)} \sum_{i,k,l=1}^n (\partial_l H_{ik}(x))^2 \, x_p \, x_q \\ 
&= \int_{\partial B_r(0)} \sum_{i,j,k,l,m=1}^n (W_{ijkl} + W_{ilkj}) \, (W_{imkl} + W_{ilkm}) \, 
x_j \, x_m \, x_p \, x_q \\ 
&= \frac{2}{n(n+2)} \, |S^{n-1}| \, \sum_{i,k,l=1}^n (W_{ipkl} + W_{ilkp}) \, (W_{iqkl} + W_{ilkq}) \, r^{n+3} \\ 
&+ \frac{1}{n(n+2)} \, |S^{n-1}| \, \sum_{i,j,k,l=1}^n (W_{ijkl} + W_{ilkj})^2 \, \delta_{pq} \, r^{n+3}. 
\end{align*}
Moreover, it follows from Corollary \ref{useful.identities} that
\begin{align*}
&\int_{\partial B_r(0)} \sum_{i,k=1}^n H_{ik}(x)^2 \, x_p \, x_q \\ 
&= \int_{\partial B_r(0)} \sum_{i,j,k,l,m,s=1}^n W_{ijkl} \, W_{imks} \, x_j \, x_l \, x_m \, x_s \, x_p \, x_q \\ 
&= \frac{2}{n(n+2)(n+4)} \, |S^{n-1}| \, \sum_{i,k,l=1}^n (W_{ipkl} + W_{ilkp}) \, (W_{iqkl} + W_{ilkq}) \, r^{n+5} \\ 
&+ \frac{1}{2n(n+2)(n+4)} \, |S^{n-1}| \, \sum_{i,j,k,l=1}^n (W_{ijkl} + W_{ilkj})^2 \, \delta_{pq} \, r^{n+5}. 
\end{align*} 
This completes the proof. \\

\begin{proposition} 
\label{integral.identity.2}
We have 
\begin{align*} 
&\int_{\partial B_r(0)} \sum_{i,k,l=1}^n (\partial_l \overline{H}_{ik}(x))^2 \, x_p \, x_q \\ 
&= \frac{2}{n(n+2)} \, |S^{n-1}| \, \sum_{i,k,l=1}^n (W_{ipkl} + W_{ilkp}) \, (W_{iqkl} + W_{ilkq}) \\ 
&\hspace{10mm} \cdot \Big [ r^{n+3} - \frac{2(n+8)}{n+4} \, r^{n+5} + \frac{n+16}{n+4} \, r^{n+7} \Big ] \\ 
&+ \frac{1}{n(n+2)} \, |S^{n-1}| \, \sum_{i,j,k,l=1}^n (W_{ijkl} + W_{ilkj})^2 \, \delta_{pq} \\ 
&\hspace{10mm} \cdot \Big [ r^{n+3} - \frac{2(n+6)}{n+4} \, r^{n+5} + \frac{n+10}{n+4} \, r^{n+7} \Big ]. 
\end{align*}
\end{proposition}

\textbf{Proof.} 
Using the identity 
\[\partial_l \overline{H}_{ik}(x) = (1 - |x|^2) \, \partial_l H_{ik}(x) - 2 \, H_{ik}(x) \, x_l\] 
and Euler's theorem, we obtain 
\begin{align*} 
&\sum_{i,k,l=1}^n (\partial_l \overline{H}_{ik}(x))^2 \\ 
&= (1 - |x|^2)^2 \, \sum_{i,k,l=1}^n (\partial_l H_{ik}(x))^2 \\ 
&- 4 \, (1 - |x|^2) \, \sum_{i,k,l=1}^n H_{ik}(x) \, x_l \, \partial_l H_{ik}(x) 
+ 4 \, |x|^2 \, \sum_{i,k=1}^n H_{ik}(x)^2 \\ 
&= (1 - |x|^2)^2 \, \sum_{i,k,l=1}^n (\partial_l H_{ik}(x))^2 - 4 \, (2 - 3 \, |x|^2) \, \sum_{i,k=1}^n H_{ik}(x)^2. 
\end{align*} 
Hence, the assertion follows from the previous proposition. \\

\begin{corollary} 
\label{integral.identity.3}
We have 
\begin{align*} 
\int_{\partial B_r(0)} \sum_{i,k,l=1}^n (\partial_l \overline{H}_{ik}(x))^2 &= \frac{1}{n} \, |S^{n-1}| \, \sum_{i,j,k,l=1}^n (W_{ijkl} + W_{ilkj})^2 \\ 
&\hspace{10mm} \cdot \Big [ r^{n+1} - \frac{2(n+4)}{n+2} \, r^{n+3} + \frac{n+8}{n+2} \, r^{n+5} \Big ]. 
\end{align*}
\end{corollary}

\vspace{2mm}

\begin{proposition} 
\label{formula.for.F}
We have 
\begin{align*} 
F(0,\varepsilon) &= -\frac{(n-2)(n+4)}{16n(n-1)(n+2)} \, |S^{n-1}| \, \sum_{i,j,k,l=1}^n (W_{ijkl} + W_{ilkj})^2 \\ 
&\hspace{10mm} \cdot \Big [ \frac{n-8}{n+4} \, \varepsilon^4 - 2 \, \varepsilon^6 + \frac{n+8}{n-10} \, \varepsilon^8 \Big ] \, \int_0^\infty (1+r^2)^{2-n} \, r^{n+3} \, dr. 
\end{align*} 
\end{proposition}

\textbf{Proof.} 
Note that $z_{(0,\varepsilon)}(x) = 0$ for all $x \in \mathbb{R}^n$. This implies 
\[F(0,\varepsilon) = -\int_{\mathbb{R}^n} \frac{n-2}{16(n-1)} \, \varepsilon^{n-2} \, (\varepsilon^2 + |x|^2)^{2-n} \, \sum_{i,k,l=1}^n (\partial_l \overline{H}_{ik}(x))^2.\] 
Using Corollary \ref{integral.identity.3}, we obtain 
\begin{align*} 
&\int_{\mathbb{R}^n} \varepsilon^{n-2} \, (\varepsilon^2 + |x|^2)^{2-n} \, \sum_{i,k,l=1}^n (\partial_l \overline{H}_{ik}(x))^2 \\ 
&= \frac{1}{n} \, |S^{n-1}| \, \sum_{i,j,k,l=1}^n (W_{ijkl} + W_{ilkj})^2 \\ 
&\hspace{10mm} \cdot \int_0^\infty (1+r^2)^{2-n} 
\, \Big [ \varepsilon^4 \, r^{n+1} - \frac{2(n+4)}{n+2} \, \varepsilon^6 \, r^{n+3} + \frac{n+8}{n+2} \, \varepsilon^8 \, r^{n+5} \Big ] \, dr. 
\end{align*} 
Moreover, we have 
\[\int_0^\infty (1+r^2)^{2-n} \, r^{n+1} \, dr = \frac{n-8}{n+2} \int_0^\infty (1+r^2)^{2-n} \, r^{n+3} \, dr\] 
and 
\[\int_0^\infty (1+r^2)^{2-n} \, r^{n+5} \, dr = \frac{n+4}{n-10} \int_0^\infty (1+r^2)^{2-n} \, r^{n+3} \, dr\] 
by Proposition \ref{integration.by.parts}. From this the assertion follows. \\

\begin{corollary}
\label{critical.point.of.F}
Assume that $n \geq 52$. Moreover, suppose that $\varepsilon_* > 0$ is defined by 
\begin{equation} 
\label{eps.star}
\bigg ( 3 + \sqrt{9 - \frac{8(n+8)(n-8)}{(n+4)(n-10)}} \bigg ) \, 
\varepsilon_*^2 = \frac{2(n-8)}{n+4}. 
\end{equation}
Then $(0,\varepsilon_*)$ is a critical point of the function $F(\xi,\varepsilon)$. Moreover, we have $\frac{\partial^2}{\partial \varepsilon^2} F(0,\varepsilon_*) > 0$.
\end{corollary} 

In the next step, we show that $(0,\varepsilon_*)$ is a strict local minimum of the function $F$. To that end, we compute the Hessian of $F$ at a point $(0,\varepsilon)$. \\

\begin{proposition} 
\label{Hessian.of.F.1}
The second order partial derivatives of the function $F(\xi,\varepsilon)$ are given by 
\begin{align*} 
\frac{\partial^2}{\partial \xi_p \, \partial \xi_q} F(0,\varepsilon) 
&= \int_{\mathbb{R}^n} (n-2)^2 \, \varepsilon^{n-2} \, (\varepsilon^2 + |x|^2)^{-n} \, \sum_{l=1}^n \overline{H}_{pl}(x) \, \overline{H}_{ql}(x) \\ 
&- \int_{\mathbb{R}^n} \frac{(n-2)^2}{4} \, \varepsilon^{n-2} \, (\varepsilon^2 + |x|^2)^{-n} \, \sum_{i,k,l=1}^n (\partial_l \overline{H}_{ik}(x))^2 \, x_p \, x_q \\ 
&+ \int_{\mathbb{R}^n} \frac{(n-2)^2}{8(n-1)} \, \varepsilon^{n-2} \, (\varepsilon^2 + |x|^2)^{1-n} \, \sum_{i,k,l=1}^n (\partial_l \overline{H}_{ik}(x))^2 \, \delta_{pq}. 
\end{align*} 
\end{proposition}

\textbf{Proof.} 
Using the identity 
\begin{align*} 
&\sum_{i,k,l=1}^n \overline{H}_{il}(x) \, \overline{H}_{kl}(x) \, \partial_i u_{(\xi,\varepsilon)}(x) \, \partial_k u_{(\xi,\varepsilon)}(x) \\ 
&= (n-2)^2 \, \varepsilon^{n-2} \, (\varepsilon^2 + |x - \xi|^2)^{-n} \, \sum_{i,k,l=1}^n \overline{H}_{il}(x) \, \overline{H}_{kl}(x) \, (x_i - \xi_i) \, (x_k - \xi_k) \\ 
&= (n-2)^2 \, \varepsilon^{n-2} \, (\varepsilon^2 + |x - \xi|^2)^{-n} \, \sum_{i,k,l=1}^n \overline{H}_{il}(x) \, \overline{H}_{kl}(x) \, \xi_i \, \xi_k, 
\end{align*} 
we obtain 
\begin{align*} 
&\frac{\partial^2}{\partial \xi_p \, \partial \xi_q} \bigg ( \frac{1}{2} \sum_{i,k,l=1}^n \overline{H}_{il}(x) \, \overline{H}_{kl}(x) \, \partial_i u_{(\xi,\varepsilon)}(x) \, \partial_k u_{(\xi,\varepsilon)}(x) \bigg ) \bigg |_{\xi=0} \\ 
&= (n-2)^2 \, \varepsilon^{n-2} \, (\varepsilon^2 + |x|^2)^{-n} \, \sum_{l=1}^n \overline{H}_{pl}(x) \, \overline{H}_{ql}(x). 
\end{align*} 
Moreover, we have 
\begin{align*} 
&\frac{\partial^2}{\partial \xi_p \, \partial \xi_q} \bigg ( \frac{n-2}{16(n-1)} \, \sum_{i,k,l=1}^n (\partial_l \overline{H}_{ik}(x))^2 \, u_{(\xi,\varepsilon)}(x)^2 \bigg ) \bigg |_{\xi=0} \\ 
&= \frac{(n-2)^2}{4} \, \varepsilon^{n-2} \, (\varepsilon^2 + |x|^2)^{-n} \, \sum_{i,k,l=1}^n (\partial_l \overline{H}_{ik}(x))^2 \, x_p \, x_q \\ 
&- \frac{(n-2)^2}{8(n-1)} \, \varepsilon^{n-2} \, (\varepsilon^2 + |x|^2)^{1-n} \, \sum_{i,k,l=1}^n (\partial_l \overline{H}_{ik}(x))^2 \, \delta_{pq}. 
\end{align*}
Finally, we have 
\begin{align*} 
&\sum_{i,k=1}^n \overline{H}_{ik}(x) \, \partial_i \partial_k u_{(\xi,\varepsilon)}(x) \\ 
&= n(n-2) \, \varepsilon^{\frac{n-2}{2}} \, (\varepsilon^2 + |x - \xi|^2)^{-\frac{n+2}{2}} \, \sum_{i,k=1}^n \overline{H}_{ik}(x) \, (x_i - \xi_i) \, (x_k - \xi_k) \\ 
&= n(n-2) \, \varepsilon^{\frac{n-2}{2}} \, (\varepsilon^2 + |x - \xi|^2)^{-\frac{n+2}{2}} \, \sum_{i,k=1}^n \overline{H}_{ik}(x) \, \xi_i \, \xi_k 
\end{align*} 
since $\overline{H}_{ik}(x)$ is trace-free. Thus, we conclude that 
\begin{align*} 
&\frac{\partial^2}{\partial \xi_p \, \partial \xi_q} \bigg ( \sum_{i,k=1}^n \overline{H}_{ik}(x) \, \partial_i \partial_k u_{(\xi,\varepsilon)}(x) \, z_{(\xi,\varepsilon)}(x) \bigg ) \bigg |_{\xi=0} \\ 
&= 2n(n-2) \, \varepsilon^{\frac{n-2}{2}} \, (\varepsilon^2 + |x|^2)^{-\frac{n+2}{2}} \, \sum_{i,k=1}^n \overline{H}_{pq}(x) \, z_{(0,\varepsilon)}(x) = 0. 
\end{align*} 
From this the assertion follows. \\

\begin{proposition} 
\label{Hessian.of.F.2}
The second order partial derivatives of the function $F(\xi,\varepsilon)$ are given by 
\begin{align*} 
&\frac{\partial^2}{\partial \xi_p \, \partial \xi_q} F(0,\varepsilon) \\ 
&= \frac{4 (n-2)^2}{n(n+2)(n+4)} \, |S^{n-1}| \, \sum_{i,k,l=1}^n (W_{ipkl} + W_{ilkp}) \, (W_{iqkl} + W_{ilkq}) \\ 
&\hspace{10mm} \cdot \Big [ \varepsilon^4 - \frac{3(n+6)}{2(n-8)} \, \varepsilon^6 \Big ] \, \int_0^\infty (1 + r^2)^{-n} \, r^{n+5} \, dr \\ 
&+ \frac{(n-2)^2}{n(n+2)(n+4)} \, |S^{n-1}| \, \sum_{i,j,k,l=1}^n (W_{ijkl} + W_{ilkj})^2 \, \delta_{pq} \\ 
&\hspace{10mm} \cdot \Big [ \varepsilon^4 - \frac{n+7}{n-8} \, \varepsilon^6 \Big ] \, \int_0^\infty (1 + r^2)^{-n} \, r^{n+5} \, dr. 
\end{align*} 
\end{proposition}

\textbf{Proof.} 
Using the identity
\begin{align*} 
&\int_{\partial B_r(0)} \sum_{l=1}^n \overline{H}_{pl}(x) \, \overline{H}_{ql}(x) \\ 
&= \int_{\partial B_r(0)} \sum_{i,j,k,l,m=1}^n W_{ipkl} \, W_{jqml} \, x_i \, x_j \, x_k \, x_m \, (1 - |x|^2)^2 \\ 
&= \frac{1}{n(n+2)} \, |S^{n-1}| \\ 
&\hspace{10mm} \cdot \sum_{i,j,k,l,m=1}^n W_{ipkl} \, W_{jqml} \, (\delta_{ij} \, \delta_{km} + \delta_{ik} \, \delta_{jm} + \delta_{im} \, \delta_{jk}) \, r^{n+3} \, (1 - r^2)^2 \\ 
&= \frac{1}{2n(n+2)} \, |S^{n-1}| \, \sum_{i,k,l=1}^n (W_{ipkl} + W_{ilkp}) \, (W_{iqkl} + W_{ilkq}) \, r^{n+3} \, (1 - r^2)^2, 
\end{align*} 
we obtain 
\begin{align*} 
&\int_{\mathbb{R}^n} \varepsilon^{n-2} 
\, (\varepsilon^2 + |x|^2)^{-n} \, \sum_{i,k,l=1}^n \overline{H}_{pl}(x) 
\, \overline{H}_{ql}(x) \\ 
&= \frac{1}{2n(n+2)} \, |S^{n-1}| \, 
\sum_{i,k,l=1}^n (W_{ipkl} + W_{ilkp}) \, (W_{iqkl} + W_{ilkq}) \\ 
&\hspace{10mm} \cdot \int_0^\infty 
(1 + r^2)^{-n} \, \Big [ \varepsilon^2 \, r^{n+3} - 2 \, \varepsilon^4 \, r^{n+5} + \varepsilon^6 \, r^{n+7} \Big ] \, dr. 
\end{align*} 
Similarly, it follows from Proposition \ref{integral.identity.2} that 
\begin{align*} 
&\int_{\mathbb{R}^n} \varepsilon^{n-2} \, (\varepsilon^2 + |x|^2)^{-n} \, \sum_{i,k,l=1}^n (\partial_l \overline{H}_{ik}(x))^2 \, x_p \, x_q \\ 
&= \frac{2}{n(n+2)} \, |S^{n-1}| \, \sum_{i,k,l=1}^n (W_{ipkl} + W_{ilkp}) \, (W_{iqkl} + W_{ilkq}) \\ 
&\hspace{10mm} \cdot \int_0^\infty (1 + r^2)^{-n} \, \Big [ \varepsilon^2 \, r^{n+3} - \frac{2(n+8)}{n+4} \, \varepsilon^4 \, r^{n+5} + \frac{n+16}{n+4} \, \varepsilon^6 \, r^{n+7} \Big ] \, dr \\ 
&+ \frac{1}{n(n+2)} \, |S^{n-1}| \, \sum_{i,j,k,l=1}^n (W_{ijkl} + W_{ilkj})^2 \, \delta_{pq} \\ 
&\hspace{10mm} \cdot \int_0^\infty (1 + r^2)^{-n} \, \Big [ \varepsilon^2 \, r^{n+3} - \frac{2(n+6)}{n+4} \, \varepsilon^4 \, r^{n+5} + \frac{n+10}{n+4} \, \varepsilon^6 \, r^{n+7} \Big ] \, dr. 
\end{align*}
Moreover, we have 
\begin{align*} 
&\int_{\mathbb{R}^n} \varepsilon^{n-2} \, (\varepsilon^2 + |x|^2)^{1-n} \, \sum_{i,k,l=1}^n (\partial_l \overline{H}_{ik}(x))^2 \, \delta_{pq} \\ 
&= \frac{1}{n} \, |S^{n-1}| \, \sum_{i,j,k,l=1}^n (W_{ijkl} + W_{ilkj})^2 \, \delta_{pq} \\ 
&\hspace{10mm} \cdot \int_0^\infty (1 + r^2)^{1-n} \, \Big [ \varepsilon^2 \, r^{n+1} - \frac{2(n+4)}{n+2} \, \varepsilon^4 \, r^{n+3} + \frac{n+8}{n+2} \, \varepsilon^6 \, r^{n+5} \Big ] \, dr. 
\end{align*} 
by Corollary \ref{integral.identity.3}. Using Proposition \ref{Hessian.of.F.1} and the identity 
\[\int_0^\infty (1 + r^2)^{1-n} \, r^{n+1} \, dr = \frac{2(n-1)}{n+2} \int_0^\infty (1 + r^2)^{-n} \, r^{n+3} \, dr,\] 
we obtain 
\begin{align*} 
&\frac{\partial^2}{\partial \xi_p \, \partial \xi_q} F(0,\varepsilon) \\ 
&= \frac{4 (n-2)^2}{n(n+2)(n+4)} \, |S^{n-1}| \, 
\sum_{i,k,l=1}^n (W_{ipkl} + W_{ilkp}) \, (W_{iqkl} + W_{ilkq}) \\ 
&\hspace{10mm} \cdot \int_0^\infty 
(1 + r^2)^{-n} \, \Big [ \varepsilon^4 \, r^{n+5} - \frac{3}{2} \, \varepsilon^6 \, r^{n+7} \Big ] \, dr \\ 
&+ \frac{(n-2)^2}{4n(n+2)} \, |S^{n-1}| \, \sum_{i,j,k,l=1}^n (W_{ijkl} + W_{ilkj})^2 \, \delta_{pq} \\ 
&\hspace{10mm} \cdot \int_0^\infty (1 + r^2)^{-n} \, \Big [ \frac{2(n+6)}{n+4} \, \varepsilon^4 \, r^{n+5} - \frac{n+10}{n+4} \, \varepsilon^6 \, r^{n+7} \Big ] \, dr \\ 
&- \frac{(n-2)^2}{8n(n-1)} \, |S^{n-1}| \, \sum_{i,j,k,l=1}^n (W_{ijkl} + W_{ilkj})^2 \, \delta_{pq} \\ 
&\hspace{10mm} \cdot \int_0^\infty (1 + r^2)^{1-n} \, \Big [ \frac{2(n+4)}{n+2} \, \varepsilon^4 \, r^{n+3} - \frac{n+8}{n+2} \, \varepsilon^6 \, r^{n+5} \Big ] \, dr.
\end{align*}
Hence, the assertion follows from the identities 
\begin{align*} 
&\int_0^\infty (1 + r^2)^{-n} \, r^{n+7} \, dr = \frac{n+6}{n-8} \int_0^\infty (1 + r^2)^{-n} \, r^{n+5} \, dr \\ 
&\int_0^\infty (1 + r^2)^{1-n} \, r^{n+3} \, dr = \frac{2(n-1)}{n+4} \int_0^\infty (1 + r^2)^{-n} \, r^{n+5} \, dr \\ 
&\int_0^\infty (1 + r^2)^{1-n} \, r^{n+5} \, dr = \frac{2(n-1)}{n-8} \int_0^\infty (1 + r^2)^{-n} \, r^{n+5} \, dr. 
\end{align*}

\vspace{2mm}

\begin{corollary} 
\label{strict.local.minimum}
Assume that $n \geq 52$ and $\varepsilon_* > 0$ is defined by (\ref{eps.star}). Then the function $F(\xi,\varepsilon)$ has a strict local minimum at the point $(0,\varepsilon_*)$.
\end{corollary}

\textbf{Proof.} 
It follows from Corollary \ref{critical.point.of.F} that $(0,\varepsilon_*)$ is a critical point of the function $F(\xi,\varepsilon)$. Moreover, we have 
$\frac{\partial^2}{\partial \varepsilon^2} F(0,\varepsilon_*) > 0$. 
Since $n \geq 52$, we have 
\[\frac{6}{n+4} < \sqrt{9 - \frac{8(n+8)(n-8)}{(n+4)(n-10)}}.\] 
This implies 
\[\frac{3(n+6)}{n+4} \, \varepsilon_*^2 < \bigg ( 3 + \sqrt{9 - \frac{8(n+8)(n-8)}{(n+4)(n-10)}} \bigg ) \, \varepsilon_*^2 = \frac{2(n-8)}{n+4}.\] 
Thus, we conclude that 
\[\frac{n+7}{n-8} \, \varepsilon_*^2 < \frac{3(n+6)}{2(n-8)} \, \varepsilon_*^2 < 1.\] 
Hence, it follows from Proposition \ref{Hessian.of.F.2} that the matrix 
$\frac{\partial^2}{\partial \xi_p \, \partial \xi_q} 
F(0,\varepsilon_*)$ is positive definite. This proves the assertion. \\

\section{Proof of the main theorem}

\begin{proposition}
\label{perturbation.argument}
Assume that $n \geq 52$. Moreover, let $g$ be a smooth metric on $\mathbb{R}^n$ of the form 
$g(x) = \exp(h(x))$, where $h(x)$ is a trace-free symmetric two-tensor on $\mathbb{R}^n$ such 
that $|h(x)| + |\partial h(x)| + |\partial^2 h(x)| \leq \alpha \leq \alpha_1$ for all $x \in \mathbb{R}^n$, 
$h(x) = 0$ for $|x| \geq 1$, and $h_{ik}(x) = \mu \, (\lambda^2 - |x|^2) \,
H_{ik}(x)$ for $|x| \leq \rho$. As above, we assume that $\lambda \leq \rho \leq 1$ and $\mu \leq 1$. If $\alpha$ and $\rho^{2-n} \, \mu^{-2} \, \lambda^{n-10}$ are 
sufficiently small, then there exists a positive function $v$ such that 
\[\Delta_g v - \frac{n-2}{4(n-1)} \, R_g \, v + n(n-2) \, v^{\frac{n+2}{n-2}} = 0,\] 
\[\int_{\mathbb{R}^n} v^{\frac{2n}{n-2}} < \Big ( \frac{Y(S^n)}{4n(n-1)} \Big )^{\frac{n}{2}},\] and  
$\sup_{|x| \leq \lambda} v(x) \geq c \, \lambda^{\frac{2-n}{2}}$. Here, $c$ is a positive constant that depends only on $n$.
\end{proposition}

\textbf{Proof.} 
By Corollary \ref{strict.local.minimum}, the function $F(\xi,\varepsilon)$ has a strict local minimum at $(0,\varepsilon_*)$. Hence, we can find an open set $\Omega' \subset \Omega$ such that $(0,\varepsilon_*) \in \Omega'$ and 
\[F(0,\varepsilon_*) < \inf_{(\xi,\varepsilon) \in \partial \Omega'} F(\xi,\varepsilon) < 0.\] 
Using Corollary \ref{key.estimate}, we obtain 
\begin{align*} 
&|\mathcal{F}_g(\lambda\xi,\lambda\varepsilon) - \lambda^8 \, \mu^2 \, F(\xi,\varepsilon)| \\ 
&\leq C \, \lambda^{\frac{8n}{n-2}} \, \mu^{\frac{2n}{n-2}} + C \, \lambda^4 \, \mu \, \Big ( \frac{\lambda}{\rho} \Big )^{\frac{n-2}{2}} + C \, \Big ( \frac{\lambda}{\rho} \Big )^{n-2} 
\end{align*} 
for all $(\xi,\varepsilon) \in \Omega$. This implies 
\begin{align*} 
&|\lambda^{-8} \, \mu^{-2} \, \mathcal{F}_g(\lambda\xi,\lambda\varepsilon) - F(\xi,\varepsilon)| \\ 
&\leq C \, \lambda^{\frac{16}{n-2}} \, \mu^{\frac{4}{n-2}} + C \, \rho^{\frac{2-n}{2}} \, \mu^{-1} \, \lambda^{\frac{n-10}{2}} 
+ C \, \rho^{2-n} \, \mu^{-2} \, \lambda^{n-10} 
\end{align*} 
for all $(\xi,\varepsilon) \in \Omega$. Hence, if $\rho^{2-n} \, \mu^{-2} \, \lambda^{n-10}$ is sufficiently small, then we have 
\[\mathcal{F}_g(0,\lambda\varepsilon_*) < \inf_{(\xi,\varepsilon) \in \partial \Omega'} \mathcal{F}_g(\lambda\xi,\lambda\varepsilon) < 0.\] 
Consequently, there exists a point $(\bar{\xi},\bar{\varepsilon}) \in \Omega'$ such that 
\[\mathcal{F}_g(\lambda\bar{\xi},\lambda\bar{\varepsilon}) = \inf_{(\xi,\varepsilon) \in \Omega'} \mathcal{F}_g(\lambda\xi,\lambda\varepsilon) < 0.\] 
By Proposition \ref{reduction.to.a.finite.dimensional.problem}, the function $v = v_{(\lambda\bar{\xi},\lambda\bar{\varepsilon})}$ is a non-negative weak solution of the partial differential equation 
\[\Delta_g v - \frac{n-2}{4(n-1)} \, R_g \, v + n(n-2) \, v^{\frac{n+2}{n-2}} = 0.\] 
Using a result of N.~Trudinger, we conclude that $v$ is smooth (see \cite{Trudinger}, Theorem 3 on p. 271). Moreover, we have 
\begin{align*}
2(n-2) \int_{\mathbb{R}^n} v^{\frac{2n}{n-2}}
&= 2(n-2) \, \Big ( \frac{Y(S^n)}{4n(n-1)} \Big )^{\frac{n}{2}} + \mathcal{F}_g(\lambda\bar{\xi},\lambda\bar{\varepsilon}) \\
&< 2(n-2) \, \Big ( \frac{Y(S^n)}{4n(n-1)} \Big )^{\frac{n}{2}}.
\end{align*}
Finally, it follows from Proposition \ref{fixed.point.argument} that $\|v - u_{(\lambda\bar{\xi},\lambda\bar{\varepsilon})}\|_{L^{\frac{2n}{n-2}}(\mathbb{R}^n)} \leq C \, \alpha$. 
This implies 
\[|B_\lambda(0)|^{\frac{n-2}{2n}} \, \sup_{|x| \leq \lambda} v(x) \geq \|v\|_{L^{\frac{2n}{n-2}}(B_\lambda(0))} \geq \|u_{(\lambda\bar{\xi},\lambda\bar{\varepsilon})}\|_{L^{\frac{2n}{n-2}}(B_\lambda(0))} - C \, \alpha.\] 
Hence, if $\alpha$ is sufficiently small, then we obtain $\lambda^{\frac{n-2}{2}} \, \sup_{|x| \leq \lambda} v(x) \geq c$. \\

\begin{proposition}
Let $n \geq 52$. Then there exists a smooth metric $g$ on $\mathbb{R}^n$ with the following properties: 
\begin{itemize}
\item[(i)] $g_{ik}(x) = \delta_{ik}$ for $|x| \geq \frac{1}{2}$ 
\item[(ii)] $g$ is not conformally flat 
\item[(iii)] There exists a sequence of non-negative smooth functions $v_\nu$ ($\nu \in \mathbb{N}$) such that  
\[\Delta_g v_\nu - \frac{n-2}{4(n-1)} \, R_g \, v_\nu + n(n-2) \, v_\nu^{\frac{n+2}{n-2}} = 0\] 
for all $\nu \in \mathbb{N}$,
\[\int_{\mathbb{R}^n} v_\nu^{\frac{2n}{n-2}} < \Big ( \frac{Y(S^n)}{4n(n-1)} \Big )^{\frac{n}{2}}\]
for all $\nu \in \mathbb{N}$, and $\sup_{|x| \leq 1} v_\nu(x) \to \infty$ as $\nu \to \infty$.
\end{itemize}
\end{proposition}

\textbf{Proof.} 
Choose a smooth cutoff function $\eta: \mathbb{R} \to \mathbb{R}$ such that $\eta(t) = 1$ for $t \leq 1$ and $\eta(t) = 0$ for $t \geq 2$. We  define a trace-free symmetric two-tensor on $\mathbb{R}^n$ by 
\[h_{ik}(x) = \sum_{N=N_0}^\infty \eta(4N^2 \, |x - y_N|) \, 2^{-N} \, (2^{-N} - |x - y_N|^2) \, H_{ik}(x - y_N),\] 
where $y_N = (\frac{1}{N},0,\hdots,0) \in \mathbb{R}^n$. It is straightforward to verify that $h(x)$ is $C^\infty$ smooth. 

Let $\alpha$ be the constant appearing in Proposition \ref{perturbation.argument}. If $N_0$ is sufficiently large, then we have $|h(x)| + |\partial h(x)| + |\partial^2 h(x)| \leq \alpha$ for all $x \in \mathbb{R}^n$ and $h(x) = 0$ for $|x| \geq \frac{1}{2}$. Moreover, we have $h_{ik}(x) = 2^{-N} \, (2^{-N} - |x - y_N|^2) \, H_{ik}(x - y_N)$ provided that $N \geq N_0$ and $|x - y_N| \leq \frac{1}{4N^2}$. Hence, we can apply Proposition \ref{perturbation.argument} with $\lambda = 2^{-N/2}$, $\mu = 2^{-N}$, and $\rho = \frac{1}{4N^2}$. From this the assertion follows. \\

\appendix

\section{An asymptotic expansion for the scalar curvature} 

Suppose that $h(x)$ is a trace-free symmetric two-tensor defined on $\mathbb{R}^n$ satisfying $|h(x)| \leq 1$ for all $x \in \mathbb{R}^n$. We define a Riemannian metric $g$ on $\mathbb{R}^n$ by $g(x) = \exp(h(x))$. In this section, we derive an approximate formula for the scalar curvature of this metric. A similar formula is derived in \cite{Ambrosetti-Malchiodi}. \\

\begin{proposition}
\label{Taylor.expansion.of.scalar.curvature}
Let $R_g$ be the scalar curvature of $g$. There exists a constant $C$, depending only on $n$, such that
\begin{align*}
&\Big | R_g - \partial_i \partial_k h_{ik} 
+ \partial_i(h_{il} \, \partial_k h_{kl}) - \frac{1}{2} \, \partial_i h_{il} \, \partial_k h_{kl} + \frac{1}{4} \, \partial_l h_{ik} \, \partial_l h_{ik} \Big | \\ 
&\leq C \, |h|^2 \, |\partial^2 h| + C \, |h| \, |\partial h|^2.
\end{align*}
\end{proposition}

\textbf{Proof.}
The Riemann curvature tensor is defined as 
\[\partial_i \Gamma_{jk}^m - \partial_j \Gamma_{ik}^m + \Gamma_{jk}^l \, \Gamma_{il}^m - \Gamma_{ik}^l \, \Gamma_{jl}^m.\] 
Hence, the scalar curvature of $g$ is given by
\[R_g = g^{jk} \, (\partial_i \Gamma_{jk}^i - \partial_j \Gamma_{ik}^i + \Gamma_{jk}^l \, \Gamma_{il}^i - \Gamma_{ik}^l \, \Gamma_{jl}^i).\] 
Since $h$ is trace-free, we have $\det g(x) = 1$ for all $x \in \mathbb{R}^n$. This implies $\Gamma_{ik}^i = \frac{1}{2} \, g^{il} \, \partial_k g_{il} = \frac{1}{2} \, \partial_k \log \det g= 0$. Therefore, we obtain 
\begin{align*} 
R_g &= g^{jk} \, \partial_i \Gamma_{jk}^i - g^{jk} \, \Gamma_{ik}^l \, \Gamma_{jl}^i \\ 
&= \partial_i(g^{jk} \, \Gamma_{jk}^i) + g^{jk} \, \Gamma_{ik}^l \, \Gamma_{jl}^i. 
\end{align*} 
Note that
\[g^{jk} \, \Gamma_{jk}^i = g^{il} \, g^{jk} \, \partial_k g_{jl}.\] 
From this it follows that 
\begin{align*}
&\Big | \partial_i(g^{jk} \, \Gamma_{jk}^i) - \partial_i \partial_k h_{ik} 
+ \frac{1}{2} \, \partial_i(h_{il} \, \partial_k h_{kl}) 
+ \frac{1}{2} \, \partial_i(h_{kl} \, \partial_k h_{il}) \Big | \\
&\leq C \, |h|^2 \, |\partial^2 h| + C \, |h| \, |\partial h|^2, 
\end{align*} 
hence 
\begin{align*}
&\Big | \partial_i(g^{jk} \, \Gamma_{jk}^i) - \partial_i \partial_k h_{ik} + \partial_i(h_{il} \, \partial_k h_{kl}) 
- \frac{1}{2} \, \partial_i h_{il} \, \partial_k h_{kl} 
+ \frac{1}{2} \, \partial_i h_{kl} \, \partial_k h_{il} \Big | \\
&\leq C \, |h|^2 \, |\partial^2 h| + C \, |h| \, |\partial h|^2, 
\end{align*} 
Moreover, we have 
\[\Big | g^{jk} \, \Gamma_{ik}^l \, \Gamma_{jl}^i 
+ \frac{1}{4} \, \partial_l h_{ik} \, \partial_l h_{ik} 
- \frac{1}{2} \, \partial_i h_{kl} \, \partial_k h_{il} \Big | \leq C \, |h| \, |\partial h|^2.\] 
Putting these facts together, we obtain
\begin{align*}
&\Big | R_g - \partial_i \partial_k h_{ik} 
+ \partial_i(h_{il} \, \partial_k h_{kl}) - \frac{1}{2} \, \partial_i h_{il} \, \partial_k h_{kl} + \frac{1}{4} \, \partial_l h_{ik} \, \partial_l h_{ik} \Big | \\ 
&\leq C \, |h|^2 \, |\partial^2 h| + C \, |h| \, |\partial h|^2.
\end{align*}
This completes the proof. \\

\section{Some useful identities}

\begin{proposition} 
\label{integration.by.parts}
Suppose that $\alpha$ and $\beta$ are real numbers satisfying $2\alpha - 2 > \beta + 1 > 0$. Then 
\[\int_0^\infty (1+r^2)^{1-\alpha} \, r^\beta \, dr = \frac{2\alpha-2}{2\alpha-\beta-3} \int_0^\infty (1+r^2)^{-\alpha} \, r^\beta \, dr\] 
and 
\[\int_0^\infty (1+r^2)^{-\alpha} \, r^{\beta+2} \, dr = \frac{\beta+1}{2\alpha-\beta-3} \int_0^\infty (1+r^2)^{-\alpha} \, r^\beta \, dr.\] 
\end{proposition}

\textbf{Proof.} Using the fundamental theorem of calculus, we obtain 
\begin{align*} 
0 &= \int_0^\infty \frac{d}{dr} \big [ (1+r^2)^{1-\alpha} \, r^{\beta+1} \big ] \, dr \\ 
&= (\beta+1) \int_0^\infty (1+r^2)^{1-\alpha} \, r^\beta \, dr - (2\alpha-2) \int_0^\infty (1+r^2)^{-\alpha} \, r^{\beta+2} \, dr. 
\end{align*}
From this the assertion follows. \\

\begin{proposition} 
Suppose that $p(x)$ is a homogenous polynomial of degree $d$. Then 
\[\int_{\partial B_1(0)} p(x) = \frac{1}{d(n+d-2)} \int_{\partial B_1(0)} \Delta p(x).\] 
\end{proposition}

\textbf{Proof.} Using the divergence theorem, we obtain 
\begin{align*} 
\int_{\partial B_1(0)} \Delta p(x) 
&= (n+d-2) \int_{B_1(0)} \Delta p(x) \\ 
&= (n+d-2) \int_{\partial B_1(0)} \sum_{k=1}^n x_k \, \partial_k p(x) \\ 
&= d(n+d-2) \int_{\partial B_1(0)} p(x). 
\end{align*}

\begin{corollary}
\label{useful.identities}
We have 
\[\int_{\partial B_1(0)} x_i \, x_j = \frac{1}{n} \, |S^{n-1}| \, \delta_{ij},\] 
\[\int_{\partial B_1(0)} x_i \, x_j \, x_k \, x_l = \frac{1}{n(n+2)} \, |S^{n-1}| \, (\delta_{ij} \, \delta_{kl} + \delta_{ik} \, \delta_{jl} + \delta_{il} \, \delta_{jk}),\] 
and 
\begin{align*} 
&\int_{\partial B_1(0)} x_i \, x_j \, x_k \, x_l \, x_p \, x_q \\ 
&= \frac{1}{n(n+2)(n+4)} \, |S^{n-1}| \, 
(\delta_{ij} \, \delta_{kl} \, \delta_{pq} + \delta_{ij} \, \delta_{kp} \, \delta_{lq} + \delta_{ij} \, \delta_{kq} \, \delta_{lp} \\ 
&\hspace{45mm} + \delta_{ik} \, \delta_{jl} \, \delta_{pq} + \delta_{ik} \, \delta_{jp} \, \delta_{lq} + \delta_{ik} \, \delta_{jq} \, \delta_{lp} \\ 
&\hspace{45mm} + \delta_{il} \, \delta_{jk} \, \delta_{pq} + \delta_{il} \, \delta_{jp} \, \delta_{kq} + \delta_{il} \, \delta_{jq} \, \delta_{kp} \\ 
&\hspace{45mm} + \delta_{ip} \, \delta_{jk} \, \delta_{lq} + \delta_{ip} \, \delta_{jl} \, \delta_{kq} + \delta_{ip} \, \delta_{jq} \, \delta_{kl} \\ 
&\hspace{45mm} + \delta_{iq} \, \delta_{jk} \, \delta_{lp} + \delta_{iq} \, \delta_{jl} \, \delta_{kp} + \delta_{iq} \, \delta_{jp} \, \delta_{kl}).
\end{align*} 
\end{corollary}


\begin{thebibliography}{99}
\bibitem{Ambrosetti}
A.~Ambrosetti, \textit{Multiplicity results for the Yamabe problem on $S^n$,} Proc. Natl. Acad. Sci. USA 99 (2002), 15252--15256

\bibitem{Ambrosetti-Malchiodi}
A.~Ambrosetti and A.~Malchiodi, \textit{A multiplicity result 
for the Yamabe problem on $S^n$,} J. Funct. Anal. 168, 529--561 (1999)

\bibitem{Aubin1} 
T.~Aubin, \textit{\'Equations diff\'erentielles non lin\'eaires et probl\`eme de Yamabe 
concernant la courbure scalaire,} J. Math. Pures Appl. 55, 269--296 (1976)

\bibitem{Aubin2} 
T.~Aubin, \textit{Sur quelques probl\`emes de courbure scalaire,} J. Funct. Anal. 240, 269--289 (2006)

\bibitem{Aubin3} 
T.~Aubin, \textit{Solution compl\`ete de la $C^0$ compacit\'e de l'ensemble des solutions de l'\'equation de Yamabe,} J. Funct. Anal. 244, 579--589 (2007)

\bibitem{Berti-Malchiodi}
M.~Berti and A.~Malchiodi, \textit{Non-compactness and multiplicity results for the Yamabe problem on $S^n$,} J. Funct. Anal. 180, 210--241 (2001)

\bibitem{Druet} 
O.~Druet, \textit{Compactness for Yamabe metrics in low dimensions,} Internat. Math. Res. Notices 23, 1143--1191 (2004)

\bibitem{Druet-Hebey1}
O.~Druet and E.~Hebey, \textit{Blow-up examples for second order elliptic PDEs of critical Sobolev growth,} Trans. Amer. Math. Soc. 357, 1915--1929 (2004)

\bibitem{Druet-Hebey2}
O.~Druet and E.~Hebey, \textit{Elliptic equations of Yamabe type,} International Mathematics Research Surveys 1, 1--113 (2005)

\bibitem{Khuri-Marques-Schoen}
M.~Khuri, F.~Marques, and R.~Schoen, \textit{A compactness theorem for the Yamabe problem,} preprint (2007)

\bibitem{Li-Zhang}
Y.Y.~Li and L.~Zhang, \textit{Compactness of solutions to the Yamabe problem II,} Calc. Var. PDE 24, 185--237 (2005)

\bibitem{Li-Zhu}
Y.Y.~Li and L.~Zhu, \textit{Yamabe type equations on three-dimensional Riemannian manifolds,} Commun. Contemp. Math. 1, 1--50 (1999)

\bibitem{Marques}
F.C.~Marques, \textit{A-priori estimates for the Yamabe problem in 
the non-locally conformally flat case,} J. Diff. Geom. 71, 315--346 (2005)

\bibitem{Pollack}
D.~Pollack, \textit{Nonuniqueness and high energy solutions for a conformally invariant scalar 
equation,} Comm. Anal. Geom. 1, 347--414 (1993)

\bibitem{Rey} 
O.~Rey, \textit{The role of the Green's function in a non-linear elliptic equation involving the critical Sobolev exponent,} J. Funct. Anal. 89, 1--52 (1990)

\bibitem{Schoen1} 
R.M.~Schoen, \textit{Conformal deformation of a Riemannian metric to constant scalar curvature,} J. Diff. Geom. 20, 479--495 (1984)

\bibitem{Schoen2} 
R.M.~Schoen, \textit{Variational theory for the total scalar curvature functional for Riemannian metrics and related topics,} Topics in the calculus of variations (ed. by Mariano Giaquinta), Lecture Notes in Mathematics, vol. 1365, 
Springer Verlag 1989, 120--154

\bibitem{Schoen3} 
R.M.~Schoen, \textit{On the number of constant scalar curvature metrics in a conformal class,} Differential geometry (ed. by H. Blaine Lawson, Jr., and Keti Tenenblat), Pitman Monographs and Surveys in Pure and Applied Mathematics, vol. 52, Longman Scientific \& Technical 1991, 311--320

\bibitem{Schoen4}
R.M.~Schoen, \textit{A report on some recent progress on nonlinear problems in geometry,} In: Surveys in differential geometry, Lehigh University, Bethlehem, PA, 1991, 201--241

\bibitem{Trudinger}
N.~Trudinger, \textit{Remarks concerning the conformal deformation of Riemannian 
structures on compact manifolds,} Annali Scuola Norm. Sup. Pisa 22, 265--274 (1968)
\end{thebibliography}
\end{document}